\documentclass[twoside,12pt]{article}
\usepackage{amssymb}

\usepackage{a4wide}
\usepackage{amsthm}
\usepackage{amsmath}
\usepackage[all]{xy}
\usepackage{enumerate}
\usepackage{fancyhdr}

\newtheorem{theorem}{Theorem}[section]
\newtheorem{proposition}[theorem]{Proposition}
\newtheorem{corollary}[theorem]{Corollary}
\newtheorem{lemma}[theorem]{Lemma}
\theoremstyle{definition}

\newtheorem*{Beweis}{Proof}
\newtheorem{definition}[theorem]{Definition}
\newtheorem{punto}[theorem]{}
\theoremstyle{remark}
\newtheorem{remark}[theorem]{Remark}

\CompileMatrices  \swapnumbers
\CompileMatrices
\input{tcilatex}

\begin{document}

\title{Duality Theorems for Crossed Products over Rings\thanks{%
MSC (2000): 16W30, 16W20, 16S35, 16S40, 16D90 \newline
Keywords: Hopf Algebras, Crossed Product, Smash Product, Duality Theorems}}
\author{\textbf{Jawad Y. Abuhlail}\thanks{%
supported by KFUPM} \\
Department of Mathematical Sciences\\
King Fahd University of Petroleum $\&$ Minerals\\
31261 Dhahran - Saudi Arabia\\
abuhlail@kfupm.edu.sa}
\date{}
\maketitle

\begin{abstract}
In this note we improve and extend duality theorems for crossed products
obtained by M. Koppinen (C. Chen) from the case of base fields (Dedekind
domains) to the case of an arbitrary Noetherian commutative ground rings
under fairly weak conditions. In particular we extend an improved version of
the celebrated Blattner-Montgomery duality theorem to the case of arbitrary
Noetherian ground rings.
\end{abstract}

\section*{Introduction}

\emph{Crossed products }in the theory of Hopf Algebras were presented
independently by R. Blattner, M. Cohen, S. Montgomery \cite{BCM86} and Y.
Doi, M. Takeuchi \cite{DT86}. The so called duality theorems for crossed
products have their roots in the theory of group rings\ (e.g. \emph{%
Cohen-Montgomery duality theorems} \cite{CM84}).

In \cite{BM85} R. Blattner and S. Montgomery extended Cohen-Montgomery
duality theorems to the case of a Hopf $R$-algebra with bijective antipode
acting on an $R$-algebra, where $R$ is a base field, providing an infinite
version of the finite one achieved independently by Van den Bergh \cite
{vdB84}. The celebrated \emph{Blattner-Montgomery} duality theorem was
extended by C. Chen and W. Nichols \cite{CN90} to the case of Dedekind
domains. In a joint paper with J. G\'{o}mez-Torrecillas and F. Lobillo 
\cite[Theorem 3.2]{AG-TL2001} that result was extended to the case of
arbitrary Noetherian ground rings.{\normalsize \vspace*{0.08cm}}

In the case of a Hopf $R$-algebra (with a not necessarily bijective
antipode) over a base field, M. Koppinen introduced in \cite[Theorem 4.2]
{Kop92} duality theorems for a right $H$\emph{-crossed product} $A\#_{\sigma
}H$ with invertible cocycle and a left $H$-module subalgebra $U\subseteq
H^{\ast }.$ For a Hopf $R$-algebra with bijective antipode and an $R$%
-subbialgebra $U\subseteq H^{\circ }$ \cite[Corollary 5.4]{Kop92} provided
an improved version of Blattner-Montgomery duality theorem, dropping the
assumption that $U\subseteq H^{\circ }$ is a Hopf $R$-subalgebra with
bijective antipode.

Inspired by the work of M. Koppinen, C. Chen{\normalsize \ }presented in 
\cite{Che93} duality theorems for right $H$-crossed products $A\#_{\sigma }H$
with invertible cocycle. Although his main results were formulated for
arbitrary ground rings, the main applications he gave were limited to the
case of a base field \cite[Corollary 4, Corollary 9]{Che93} or a Dedekind
domain \cite[Corollary 5, Corollary 10]{Che93}.

The main objective of this note is unify these duality theorems and their
proofs as well as to generalize them to the case of arbitrary Noetherian
ground rings under fairly weak conditions. Another improvement is weakening
the assumption that the antipode of the Hopf algebra $H$ is bijective by
replacing it with the weaker condition that $H$ has a \emph{twisted antipode,%
} i.e. $H^{op}$ has an antipode $\overline{S}.$

In the first section we present the needed definitions and Lemmata. In the
second section we present the main result (Theorem \ref{haupt}) for a Hopf $%
R $-algebra with twisted antipode, a right $H$-crossed product $A\#_{\sigma
}H$ with invertible cocycle and a right $H$-module $R$-subalgebra $%
U\subseteq H^{\ast },$ where $R$ is an arbitrary Noetherian ground ring. In
case $_{R}H$ is locally projective we introduce a right $H$-submodule $%
H^{\upsilon }\subseteq H^{\ast },$ such that $(H^{\upsilon },H)$ satisfies
the modified RL-condition (\ref{RL}) with respect to $H$ and use it to
present results analog to those of M. Koppinen \cite{Kop92} (Theorem \ref
{v-hm1} and Corollary \ref{normal}).

As a corollary, with $\sigma $ trivial, Theorem \ref{BM} generalizes
Koppinen's version of the Blattner-Montgomery duality theorem 
\cite[Corollary 5.4]{Kop92} to the case of arbitrary Noetherian ground rings
(this improves also \cite[Theorem 3.2]{AG-TL2001}). Corollary \ref{BM-dense}
extends \cite[Corollary 9.4.11]{Mon93} to the case of arbitrary QF ground
rings (for an arbitrary right $H$-crossed product see Corollary \ref{dense}%
). Given a Hopf $R$-algebra $H$ with twisted antipode and a right $H$%
-crossed product $A\#_{\sigma }H$ with invertible cocycle, Theorem \ref
{cleft-du} provides a version Theorem \ref{haupt} formulated for the \emph{%
cleft }$H$\emph{-extension }$(A\#_{\sigma }H)/A.$

The third section deals with the case of an arbitrary Hopf $R$-algebra (not
necessarily with twisted antipode). There we generalize results of C. Chen 
\cite{Che93} from the case of a base field or a Dedekind domain to the case
of arbitrary Noetherian ring. For a locally projective Hopf $R$-algebra $H,$
we consider the $R$-subalgebra $H^{\omega }\subseteq H^{\ast }$ presented by
M. Koppinen and prove his main duality theorem \cite[Theorem 4.2]{Kop92}
over arbitrary Noetherian ground rings. We also generalize several
corollaries of \cite[Section 5]{Kop92} to the case of arbitrary Noetherian
ground rings.

With $R$ we denote a commutative ring with $1_{R}\neq 0_{R}.$ The category
of unital $R$-(bi)modules will be denoted by $\mathcal{M}_{R}.$ We consider $%
R$ as a linear topological ring with the discrete topology. For $R$-modules $%
M,$ $N$ we say an $R$-submodule $K\subset M$ is $N$\emph{-pure}\textbf{,} if
the canonical map $id_{K}\otimes \iota _{N}:K\otimes _{R}N\rightarrow
M\otimes _{R}N$ is injective. If $K\subset M$ is $N$-pure for every $R$%
-module $N,$ then we say $K\subset M$ is \emph{pure}\textbf{\ }(in the sense
of Cohn). For $R$-modules $M,N$ we denote by $\tau :M\otimes
_{R}N\rightarrow N\otimes _{R}M$ the canonical \emph{twist} isomorphism. Let 
$A$ be an $R$-algebra. The category of \emph{unital }left (resp. right) $A$%
-modules will be denoted by $_{A}\mathcal{M}$ (resp. $\mathcal{M}_{A}$). For
an $A$-module $M$ we call an $A$-submodule $K\subset M$ $R$\emph{-cofinite},
if $M/K$ is finitely generated in $\mathcal{M}_{R}.$

We assume the reader is familiar with the theory of Hopf $R$-algebras. For
the basic definitions and concepts we refer to \cite{Swe69} and \cite{Mon93}%
. For an $R$-coalgebra $(C,\Delta _{C},\varepsilon _{C})$ we call a \emph{%
pure }$R$-submodule $\widetilde{C}\subseteq C$ an $R$-coalgebra provided $%
\Delta _{C}(\widetilde{C})\subseteq \widetilde{C}\otimes _{R}\widetilde{C}.$
For an $R$-coalgebra $C\;$and an $R$-algebra $A,$ we consider $\mathrm{Hom}%
_{R}(C,A)$ as an $R$-algebra under the so called \emph{convolution product }$%
(f\star g)(c):=\sum f(c_{1})g(c_{2})$ and unity $\eta _{A}\circ \varepsilon
_{C}.$

For an $R$-coalgebra $C$ and a right $C$-comodule $(M,\varrho _{M})$ we
denote by $\mathrm{Cf}(M)\subseteq C$ the $R$-submodule generated by $%
\{m_{<0>}\mid m\in M,$ $\varrho _{M}(m)=\sum m_{<0>}\otimes m_{<1>}\}.$ For
an $R$-bialgebra $H$ and a right $H$-comodule $M,$ we set $M^{coH}=\{m\in
M\mid \varrho _{M}(m)=m\otimes 1_{H}\}.$

\section{Preliminaries}

In this section we introduce the needed definitions and results.

\begin{punto}
\textbf{Measuring }$R$-\textbf{pairings. }Let $C$ be an $R$-coalgebra and $A$
be an $R$-algebra with a morphism of $R$-algebras $\beta :A\rightarrow $ $%
C^{\ast },$ $a\mapsto \lbrack c\mapsto <a,c>].$ Then we call $P:=(A,C)$ a 
\emph{measuring }$R$\emph{-pairing} (the terminology is inspired by 
\cite[Page 139]{Swe69}). In this case $C$ is an $A$-bimodule through the
left and the right $A$-actions 
\begin{equation}
a\rightharpoonup c:=\sum c_{1}<a,c_{2}>\text{ and }c\leftharpoonup a:=\sum
<a,c_{1}>c_{2}\text{ for all }a\in A,\text{ }c\in C.  \label{C-r}
\end{equation}
\end{punto}

\begin{punto}
\textbf{The }$\alpha $\textbf{-condition}. Let $V,W$ be $R$-modules with an $%
R$-linear map $\beta :V\rightarrow W^{\ast }.$ We say the $R$\emph{-pairing} 
$P:=(V,W)$ satisfies the $\alpha $\emph{-condition }(or $P$ is an $\alpha $%
\emph{-pairing}), if for every $R$-module $M$ the following map is
injective: 
\begin{equation}
\alpha _{M}^{P}:M\otimes _{R}W\ \rightarrow \mathrm{Hom}_{R}(V,M),\text{ }%
\sum m_{i}\otimes w_{i}\mapsto \lbrack v\mapsto \sum m_{i}<v,w_{i}>].
\label{alpha}
\end{equation}
We say an $R$-module $W$ satisfies the $\alpha $\emph{-condition}, if the
canonical $R$-pairing $(W^{\ast },W)$ satisfies the $\alpha $-condition
(equivalently, if $_{R}W$ is locally projective in the sense of B.
Zimmermann-Huisgen \cite[Theorem 2.1]{Z-H76}, \cite[Theorem 3.2]{Gar76}). If 
$_{R}W$ is locally projective, then $_{R}W$ is flat and $R$-cogenerated
(e.g. \cite[Bemerkung 2.1.5]{Abu2001}).
\end{punto}

\begin{punto}
\textbf{The }$C$\textbf{-adic topology.}\label{adic-top} Let $P=(A,C)$ be a
measuring $R$-pairing and consider $C$ as a left $A$-module with the induced
left $A$-action in (\ref{C-r}). Then $A$ becomes a left linear topological $R
$-algebra under the so called $C$\emph{-adic topology} $\mathcal{T}_{C-}(A)$
with neighbourhood basis of $0_{A}:$%
\begin{equation*}
\mathcal{B}_{C-}(0_{A})=\{(0_{C}:W)\mid W\subset C\text{ is a finite subset}%
\}.
\end{equation*}
The category of \emph{discrete} left $(A,\mathcal{T}_{C-}(A))$-modules is
denoted by $\sigma \lbrack _{A}C].$ In fact $\sigma \lbrack _{A}C]$ is the 
\emph{smallest }Grothendieck full subcategory of $_{A}\mathcal{M}$ that
contains $C.$ The reader is referred to \cite{AW97}, \cite{Ber94} for more
investigation of this topology and to \cite{Wis91} for the well developed
theory of categories of type $\sigma \lbrack M].$
\end{punto}

\begin{punto}
\label{rat-dar}Let $P=(A,C)$ be a measuring $\alpha $-pairing. Let $M$ be a
left $A$-module and consider the canonical $A$-linear map $\rho
_{M}:M\rightarrow \mathrm{Hom}_{R}(A,M).$ We set $\mathrm{Rat}%
^{C}(_{A}M):=\rho _{M}^{-1}(M\otimes _{R}C)$ and call $M$ $C$\emph{-rational}%
, if $\mathrm{Rat}^{C}(_{A}M)=M.$ If $_{A}M$ is $C$-rational, then we have
an $R$-linear map $\varrho _{M}:=(\alpha _{M}^{P})^{-1}\circ \rho
_{M}:M\rightarrow M\otimes _{R}C.$ The class of $C$-rational left $A$%
-modules build a \emph{full }subcategory of $_{A}\mathcal{M},$ which we
denote with $\mathrm{Rat}^{C}(_{A}\mathcal{M})$ (see \cite[Lemma 2.2.7]
{Abu2001}).
\end{punto}

\begin{theorem}
\label{equal}\emph{(\cite[Theorems 1.14, 1.15]{Abu-b})} For a measuring $R$%
-pairing $P=(${$A$}$,${$C$}$)$ the following are equivalent:

\begin{enumerate}
\item  $P$ satisfies the $\alpha $-condition;

\item  $_{R}C$ is locally projective and $\beta _{P}(A)\subseteq C^{\ast }$
is dense.
\end{enumerate}

If these equivalent conditions are satisfied, then we have isomorphisms of
categories 
\begin{equation*}
\begin{tabular}{lllll}
$\mathcal{M}^{{C}}$ & $\simeq $ & $\sigma \lbrack _{A}C]$ & $=$ & $\mathrm{%
Rat}^{{C}}(_{A}\mathcal{M})$ \\ 
& $\simeq $ & $\sigma \lbrack _{C^{\ast }}C]$ & $=$ & $\mathrm{Rat}^{{C}%
}(_{C^{\ast }}\mathcal{M}).$%
\end{tabular}
\end{equation*}
\end{theorem}

\begin{punto}
\label{du-co}(\cite[Remark 2.14, Proposition 2.15]{AG-TL2001}). Assume $R$
to be Noetherian. Let $A$ be an $R$-algebra and consider $A^{\ast }$ as an $A
$-bimodule through the regular left and right actions 
\begin{equation}
(af)(b)=f(ba)\text{ and }(fa)(b)=f(ab)\text{ for }a,b\in A\text{ and }f\in
H^{\ast }.  \label{regular}
\end{equation}
We set 
\begin{equation*}
\begin{tabular}{lll}
$A^{\circ }$ & $:=$ & $\{f\in A^{\ast }\mid AfA\text{ is finitely generated
in }\mathcal{M}_{R}\}$ \\ 
& $=$ & $\{f\in A^{\ast }\mid $ $\mathrm{Ke}(f)$ contains an $R$-cofinite $A$%
-ideal$\}.$%
\end{tabular}
\end{equation*}
Then $(A,A^{\circ })$ is a measuring $\alpha $-pairing if and only if $%
A^{\circ }\subset R^{A}$ is pure. In this case $A^{\circ }$ is a locally
projective $R$-coalgebra and for every $R$-subcoalgebra $\widetilde{C}%
\subseteq A^{\circ },$ the induced $R$-pairing $(A,\widetilde{C})$ is a
measuring $\alpha $-pairing.

An $R$-algebra (resp. an $R$-bialgebra, a Hopf $R$-algebra) $A$ with $%
A^{\circ }\subset R^{A}$ pure will be called an $\alpha $\emph{-algebra}
(resp. an $\alpha $\emph{-bialgebra}, a \emph{Hopf }$\alpha $\emph{-algebra}%
). If $H$ is an $\alpha $-bialgebra (resp. a Hopf $\alpha $-algebra), then $%
H^{\circ }$ is an $R$-bialgebra (resp. a Hopf $R$-algebra).
\end{punto}

\begin{punto}
\label{AsgH}(\cite{BCM86}, \cite{DT86}) Let $H\;$be an $R$-bialgebra and $A\;
$an $R$-algebra. A \emph{weak left }$H$\emph{-action }on $A$ is an $R$%
-linear map $w:H\otimes _{R}A\rightarrow A,$ $h\otimes a\mapsto ha,$ such
that the induced $R$-linear map $\beta :A\rightarrow \mathrm{Hom}_{R}(H,A),$ 
$a\mapsto \lbrack h\mapsto ha]$ is an $R$-algebra morphism and $%
1_{H}\rightharpoonup a=a$ for all $a\in A.$

Let $A$ have a weak left $H$-action and $\sigma :H\otimes _{R}H\rightarrow A$
an $R$-linear map. Then $A\#_{\sigma }H:=A\otimes _{R}H$ is a (\emph{not
necessarily associative}) $R$-algebra under the multiplication 
\begin{equation}
(a\#_{\sigma }h)(\widetilde{a}\#_{\sigma }\widetilde{h}):=\sum a(h_{1}%
\widetilde{a})\sigma (h_{2}\otimes \widetilde{h}_{1})\#_{\sigma }h_{3}%
\widetilde{h}_{2}  \label{cp}
\end{equation}
and has in general no unity. If $A\#_{\sigma }H$ is an associative\emph{\ }$R
$-algebra with unity $1_{A}\#_{\sigma }1_{H},$ then $A\#_{\sigma }H$ is
called a \emph{right }$H$\emph{-crossed product.}\textbf{\ }In this case $%
(A\#_{\sigma }H,id\otimes \Delta _{H})$ is a right $H$-comodule algebra with 
$(A\#_{\sigma }H)^{coH}=A.$\textbf{\ }We say $\sigma $ in invertible, if
it's invertible in $(\mathrm{Hom}_{R}(H\otimes _{R}H,A),\star ).$
\end{punto}

\begin{lemma}
\emph{(\cite{BCM86}, \cite[Lemma 10]{DT86})}\label{cp-alg} Let $H\;$be an $R$%
-bialgebra, $A$ an $R$-algebra with a weak left $H$-action and $\sigma \in 
\mathrm{Hom}_{R}(H\otimes _{R}H,A).$

\begin{enumerate}
\item  $1\#_{\sigma }1$ is a unity for $A\#_{\sigma }H$ if and only if $%
\sigma $ is \emph{normal} \emph{(}i.e. 
\begin{equation}
\sigma (h\otimes 1_{H})=\varepsilon (h)1_{A}=\sigma (1_{H}\otimes h)\text{
for all }h\in H\emph{)}.
\end{equation}

\item  Assume $\sigma $ to be normal. Then $A\#_{\sigma }H$ is an
associative $R$-algebra if and only if $\sigma $ is a \emph{cocycle} \emph{(}%
i.e. 
\begin{equation}
\sum [h_{1}\sigma (k_{1}\otimes l_{1})]\sigma (h_{2}\otimes k_{2}l_{2})=\sum
\sigma (h_{1}\otimes k_{1})\sigma (h_{2}k_{2}\otimes l)\text{ for all }%
h,k,l\in H\text{\emph{)}}  \label{cocy}
\end{equation}
and satisfies the \emph{twisted module condition} 
\begin{equation}
\sum [h_{1}[k_{1}a]]\sigma (h_{2}\otimes k_{2})=\sum \sigma (h_{1}\otimes
k_{1})[(h_{2}k_{2})a]\text{ for all }h,k\in H,\text{ }a\in A.  \label{twist}
\end{equation}
\end{enumerate}
\end{lemma}

\begin{punto}
\textbf{Left smash product}. Let $H$ be an $R$-bialgebra and $A$ a left $H$%
-module algebra. Then 
\begin{equation*}
\sigma :H\otimes _{R}H\rightarrow A,\text{ }h\otimes k\mapsto \varepsilon
(h)\varepsilon (k)1_{A}
\end{equation*}
is a \emph{trivial }normal cocycle and satisfies the twisted module
condition (\ref{twist}). By Lemma \ref{cp-alg} $A\#H:=A\#_{\sigma }H$ is an
associative $R$-algebra with multiplication 
\begin{equation}
(a\#h)\bullet (\widetilde{a}\#\widetilde{h})=\sum a(h_{1}\widetilde{a}%
)\#h_{2}\widetilde{h}  \label{s-ls}
\end{equation}
and unity $1_{A}\#1_{H}.$ If the left $H$-action on $A$ is also \emph{trivial%
}, then $A\#H=A\otimes _{R}H$ as $R$-algebras. The $R$-algebra $A\#H$ was
presented by M. Sweedler \cite[Pages 155-156]{Swe69}.
\end{punto}

\section{The main Duality Theorem}

\qquad In this section we present the main result in this note, namely
Theorem \ref{haupt}. For the convention of the reader we begin with some
definitions.

\begin{punto}
\label{smash}(\cite{Doi84}, \cite[Page 375]{Doi92})\textbf{\ }Let $H$ be an $%
R$-bialgebra and $B$ a right $H$-comodule algebra. Then $\#(H,B):=(\mathrm{%
Hom}_{R}(H,B),\widehat{\star })$ is an associative $R$-algebra with
multiplication 
\begin{equation}
(f\widehat{\star }g)(h)=\sum f(g(h_{2})_{<1>}h_{1})g(h_{2})_{<0>}\text{ for
all }f,g\in \mathrm{Hom}_{R}(H,B),\text{ }h\in H  \label{brs}
\end{equation}
and unity $\eta _{B}\circ \varepsilon _{H}.$ If $U\subseteq H^{\ast }$ is a 
\emph{right} $H$-module subalgebra (with $\varepsilon _{H}\in U$), then $%
B\#U:=B\otimes _{R}U$ is an associative $R$-algebra with multiplication 
\begin{equation}
(b\#f)(\widetilde{b}\#\widetilde{f})=\sum b\widetilde{b}_{<0>}\#(f\widetilde{%
b}_{<1>})\star \widetilde{f}\text{ for all }b,\widetilde{b}\in B,\text{ }f,%
\widetilde{f}\in U  \label{r-smash}
\end{equation}
(and unity $1_{B}\#\varepsilon _{H}$).
\end{punto}

\begin{remark}
Let $R$ be Noetherian, $H$ an $\alpha $-bialgebra, $U\subseteq H^{\circ }$
an $R$-subbialgebra and consider the $\alpha $-pairing $P:=(H,U).$ Since $H$
is a left $U$-module algebra under the action $f\rightharpoonup h:=\sum
h_{1}f(h_{2}),$ we can endow $H\otimes _{R}U$ with the structure of a \emph{%
left smash algebra} under the multiplication (\ref{s-ls}). On the other hand 
$H$ is a right $H$-comodule algebra under $\Delta _{H},$ $U\subseteq H^{\ast
}$ is a right $H$-module subalgebra under the right regular $H$-action (\ref
{regular}) and $H\otimes _{R}U$ can be endowed with the structure of a \emph{%
right smash algebra} under the multiplication (\ref{r-smash}). It can be
easily seen that the two $R$-algebras are isomorphic. In fact we have for
arbitrary $h,\widetilde{h}\in H,$ $f,\widetilde{f}\in U$ and all $k\in H:$%
\begin{equation*}
\begin{tabular}{lll}
$\alpha _{H}^{P}((h\#f)\bullet (\widetilde{h}\#\widetilde{f}))(k)$ & $=$ & $%
\alpha _{H}^{P}(\sum h(f_{1}\rightharpoonup \widetilde{h})\#f_{2}\star 
\widetilde{f})(k)$ \\ 
& $=$ & $\sum h(f_{1}\rightharpoonup \widetilde{h})(f_{2}\star \widetilde{f}%
)(k)$ \\ 
& $=$ & $\sum h\widetilde{h}_{1}f_{1}(\widetilde{h}_{2})f_{2}(k_{1})%
\widetilde{f}(k_{2})$ \\ 
& $=$ & $\sum h\widetilde{h}_{1}f(\widetilde{h}_{2}k_{1})\widetilde{f}(k_{2})
$ \\ 
& $=$ & $\sum h\widetilde{h}_{1}(f\widetilde{h}_{2})(k_{1})\widetilde{f}%
(k_{2})$ \\ 
& $=$ & $\alpha _{H}^{P}(\sum h\widetilde{h}_{1}\#(f\widetilde{h}_{2})\star 
\widetilde{f})(k)$ \\ 
& $=$ & $\alpha _{H}^{P}((h\#f)(\widetilde{h}\#\widetilde{f}))(k).$%
\end{tabular}
\end{equation*}
Since $\alpha _{H}^{P}$ is injective, we get $(h\#f)\bullet (\widetilde{h}\#%
\widetilde{f})=(h\#f)(\widetilde{h}\#\widetilde{f})$ and we are done.
\end{remark}

The following definition provides a generalization of the \emph{RL-condition}
suggested by \cite{BM85}:

\begin{definition}
Let $H$ be an $R$-bialgebra, $U\subseteq H^{\ast }$ a right $H$-module
subalgebra under the right regular $H$-action, $V\subseteq H^{\ast }$ an $R$%
-submodule and consider the $R$-linear maps 
\begin{equation}
\begin{tabular}{llllllll}
$\lambda $ & $:$ & $H\#U$ & $\rightarrow $ & $\mathrm{End}_{R}(H),$ & $\sum
h_{j}\#g_{j}$ & $\mapsto $ & $[\widetilde{k}\mapsto
h_{j}(g_{j}\rightharpoonup \widetilde{k})],$ \\ 
$\rho $ & $:$ & $V$ & $\rightarrow $ & $\mathrm{End}_{R}(H),$ & $g$ & $%
\mapsto $ & $[\widetilde{k}\mapsto \widetilde{k}\leftharpoonup g].$%
\end{tabular}
\label{lam-ro}
\end{equation}
We say $(V,U)$ satisfies the \emph{RL-condition} with respect to $H,$
provided $\rho (V)\subseteq \lambda (H\#U),$ i.e. if 
\begin{equation}
\text{for every }g\in V,\text{ }\exists \text{ }\{(h_{j},g_{j})\}\subset
H\times U,\text{ s.t. }\widetilde{k}\leftharpoonup g=\sum
h_{j}(g_{j}\rightharpoonup \widetilde{k})\text{ for all }\widetilde{k}\in H.
\label{RL}
\end{equation}
We say $U$ satisfies the RL-condition with respect to $H,$ if $(U,U)$
satisfies the RL-condition with respect to $H.$
\end{definition}

\begin{lemma}
\label{density}Let $H$ be an $R$-bialgebra, $U\subseteq H^{\ast }$ a right $H
$-module subalgebra and consider $H$ as a right $H$-comodule algebra through 
$\Delta _{H}.$ Let $\#(H,H)$ and $H\#U$ be the $R$-algebras defined in \emph{%
\ref{smash}} and consider the canonical $R$-algebra morphism $\beta
:H\#U\rightarrow \#(H,H).$

\begin{enumerate}
\item  If $_{R}H$ is finitely generated projective, then $H\#H^{\ast }%
\overset{\beta }{\simeq }\#(H,H)$ as $R$-algebras.

\item  If $H$ is a Hopf $R$-algebra with twisted antipode, then $%
\#(H,H)\simeq \mathrm{End}_{R}(H)$ as $R$-algebras.

\item  Let $H$ be a finitely generated projective Hopf $R$-algebra. Then $%
\lambda :H\#H^{\ast }\rightarrow \mathrm{End}_{R}(H),$ defined in \emph{(\ref
{lam-ro}),} is an $R$-algebra isomorphism. In particular $H^{\ast }$
satisfies the RL-condition \emph{(\ref{RL})} with respect to $H.$

\item  If $_{R}H$ is locally projective and $U\subseteq H^{\ast }$ is dense,
then\emph{\ }$\beta (H\#U)\subseteq \#(H,H)$ is a dense $R$-subalgebra. If
moreover $H$ is a Hopf $R$-algebra with twisted antipode and $_{R}H$ is
projective, then $H\#U\overset{\lambda }{\hookrightarrow }\mathrm{End}_{R}(H)
$ is a dense $R$-subalgebra.
\end{enumerate}
\end{lemma}

\begin{Beweis}
\begin{enumerate}
\item  Since $_{R}H$ is finitely generated projective, $\beta $ is bijective.

\item  Let $H$ be a Hopf $R$-algebra with twisted antipode $\overline{S}$
and consider the $R$-linear maps 
\begin{equation*}
\begin{tabular}{llllllllll}
$\phi _{1}$ & $:$ & $\#(H,H)$ & $\rightarrow $ & $\mathrm{End}_{R}(H),$ & $f$
& $\mapsto $ & $[h$ & $\mapsto $ & $\sum f(h_{2})h_{1}],$ \\ 
$\phi _{2}$ & $:$ & $\mathrm{End}_{R}(H)$ & $\rightarrow $ & $\#(H,H),$ & $g$
& $\mapsto $ & $[k$ & $\mapsto $ & $\sum g(k_{2})\overline{S}(k_{1})].$%
\end{tabular}
\end{equation*}

For arbitrary $f,g\in \#(H,H)$ and $h\in H$ we have 
\begin{equation*}
\begin{tabular}{lllll}
$\phi _{1}(f\widehat{\star }g)(h)$ & $=$ & $\sum (f\widehat{\star }%
g)(h_{2})h_{1}$ & $=$ & $\sum f(g(h_{3})_{2}h_{2})g(h_{3})_{1}h_{1}$ \\ 
& $=$ & $\sum f(g(h_{2})_{2}h_{12})g(h_{2})_{1}h_{11}$ & $=$ & $\phi
_{1}(f)(\sum g(h_{2})h_{1})$ \\ 
& $=$ & $(\phi _{1}(f)\circ \phi _{1}(g))(h),$ &  & 
\end{tabular}
\end{equation*}
i.e. $\phi _{1}$ is an $R$-algebra morphism. For all $R$-linear maps $%
f,g:H\rightarrow H$ we have 
\begin{equation*}
\begin{tabular}{lllll}
$(\phi _{1}\circ \phi _{2})(g)(h)$ & $=$ & $\sum \phi _{2}(g)(h_{2})h_{1}$ & 
$=$ & $\sum g(h_{3})\overline{S}(h_{2})h_{1}$ \\ 
& $=$ & $\sum g(h_{2})\varepsilon (h_{1})$ & $=$ & $g(h),$ \\ 
$(\phi _{2}\circ \phi _{1})(f)(h)$ & $=$ & $\sum \phi _{1}(f)(h_{2})%
\overline{S}(h_{1})$ & $=$ & $\sum f(h_{3})h_{2}\overline{S}(h_{1})$ \\ 
& $=$ & $\sum f(h_{2})\varepsilon (h_{1})$ & $=$ & $f(h).$%
\end{tabular}
\end{equation*}
\newline
Hence $\phi _{1}$ is an $R$-algebra isomorphism with inverse $\phi _{2}.$

\item  Let $H$ be a finitely generated projective Hopf $R$-algebra. By (1) $%
\&\;$(2) we have $R$-algebra isomorphisms $H\#U\overset{\beta }{\simeq }%
\#(H,H)\overset{\phi _{1}}{\simeq }\mathrm{End}_{R}(H)$ (recall that the
antipode of a finitely generated projective Hopf $R$-algebra is bijective by 
\cite[Proposition 4]{Par71}, hence $H$ has a twisted antipode $\overline{S}%
:=S^{-1}$). So $\lambda =\phi _{1}\circ \beta :H\#U\rightarrow \mathrm{End}%
_{R}(H)$ is an $R$-algebra isomorphisms. In particular $\rho (H^{\ast
})\subseteq \mathrm{End}_{R}(H)=\lambda (H\#H^{\ast }),$ i.e. $H^{\ast }$
satisfies the RL-condition (\ref{RL}) with respect to $H.$

\item  By \cite[Corollary 3.20]{Abu2002} $\beta (H\#U)\subseteq \#(H,H)$ is
a dense $R$-subalgebra. If $H$ is a Hopf $R$-algebra with twisted antipode
then $\#(H,H)\overset{\phi _{1}}{\simeq }\mathrm{End}_{R}(H)$ as $R$%
-algebras by (2) and we are done (notice that $\beta $ is an embedding, if $%
_{R}H$ is projective).$\blacksquare $
\end{enumerate}
\end{Beweis}

\begin{lemma}
\label{alp=bet}Let $H$ be a Hopf $R$-algebra with twisted antipode, $A$ an $R
$-algebra, $U\subseteq H^{\ast }$ an $R$\emph{-submodule} and consider the $R
$-pairing $P:=(H,U).$ Then the canonical $R$-linear map $\alpha :=\alpha
_{A\otimes _{R}H}^{P}:(A\otimes _{R}H)\otimes _{R}U\rightarrow \mathrm{Hom}%
_{R}(H,A\otimes _{R}H)$ is injective if and only if the following map is
injective 
\begin{equation}
\chi :A\otimes _{R}(H\otimes _{R}U)\rightarrow \mathrm{End}_{-A}(H\otimes
_{R}A),\text{ }a\otimes (h\otimes f)\mapsto \lbrack (k\otimes \widetilde{a}%
)\mapsto h(f\rightharpoonup k)\otimes a\widetilde{a}].  \label{bet}
\end{equation}
\end{lemma}

\begin{Beweis}
Assume $H$ to have a twisted antipode $\overline{S}.$ We show first that the 
$R$-linear map 
\begin{equation*}
\epsilon :\mathrm{Hom}{\normalsize _{R}}(H,A\otimes _{R}H)\rightarrow 
\mathrm{End}_{-A}(H\otimes _{R}A),\text{ }g\mapsto \lbrack k\otimes 
\widetilde{a}\mapsto \tau (g(k_{2}))(k_{1}\otimes \widetilde{a})]
\end{equation*}
is bijective with inverse 
\begin{equation*}
\epsilon ^{-1}:\mathrm{End}_{-A}(H\otimes _{R}A)\rightarrow \mathrm{Hom}%
{\normalsize _{R}}(H,A\otimes _{R}H),\text{ }f\mapsto \lbrack k\mapsto \tau
(f(k_{2}\otimes 1_{A}))(1_{A}\otimes \overline{S}(k_{1}))].
\end{equation*}
In fact we have for all $f\in \mathrm{End}_{-A}(H\otimes _{R}A),$ $k\in H,$ $%
\widetilde{a}\in A:$%
\begin{equation*}
\begin{tabular}{lll}
$\epsilon (\epsilon ^{-1}(f))(k\otimes \widetilde{a})$ & $=$ & $\sum \tau
\lbrack \epsilon ^{-1}(f)(k_{2})](k_{1}\otimes \widetilde{a})$ \\ 
& $=$ & $\sum \tau \lbrack \tau (f(k_{22}\otimes 1_{A}))(1_{A}\otimes 
\overline{S}(k_{21}))](k_{1}\otimes \widetilde{a})$ \\ 
& $=$ & $\sum f(k_{22}\otimes 1_{A})(\overline{S}(k_{21})\otimes
1_{A})(k_{1}\otimes \widetilde{a})$ \\ 
& $=$ & $\sum f(k_{2}\otimes 1_{A})(\overline{S}(k_{12})k_{11}\otimes 
\widetilde{a})$ \\ 
& $=$ & $\sum f(k_{2}\otimes 1_{A})(\varepsilon _{H}(k_{1})1_{H}\otimes 
\widetilde{a})$ \\ 
& $=$ & $f(k\otimes \widetilde{a})$%
\end{tabular}
\end{equation*}
and for all $g\in \mathrm{Hom}_{R}(H,A\otimes _{R}H),$ $k\in H:$%
\begin{equation*}
\begin{tabular}{lll}
$\epsilon ^{-1}(\epsilon (g))(k)$ & $=$ & $\sum \tau \lbrack \epsilon
(g)(k_{2}\otimes 1_{A})](1_{A}\otimes \overline{S}(k_{1}))$ \\ 
& $=$ & $\sum \tau \lbrack \tau (g(k_{22}))(k_{21}\otimes
1_{A})](1_{A}\otimes \overline{S}(k_{1}))$ \\ 
& $=$ & $\sum g(k_{22})(1_{A}\otimes k_{21}\overline{S}(k_{1}))$ \\ 
& $=$ & $\sum g(k_{2})(1_{A}\otimes k_{12}\overline{S}(k_{11}))$ \\ 
& $=$ & $\sum g(k_{2})(1_{A}\otimes \varepsilon _{H}(k_{1}))$ \\ 
& $=$ & $g(k).$%
\end{tabular}
\end{equation*}

Moreover we have for all $a\in A,$ $h\in H,$ $f\in U$ and $k\in H:$%
\begin{equation*}
\begin{tabular}{lll}
$(\epsilon \circ \alpha )(a\otimes (h\otimes f))(k\otimes \widetilde{a})$ & $%
=$ & $\tau (\alpha (a\otimes (h\otimes f))(k_{2}))(k_{1}\otimes \widetilde{a}%
)$ \\ 
& $=$ & $\sum f(k_{2})(h\otimes a)(k_{1}\otimes \widetilde{a})$ \\ 
& $=$ & $\sum hf(k_{2})k_{1}\otimes a\widetilde{a}$ \\ 
& $=$ & $h(f\rightharpoonup k)\otimes a\widetilde{a}$ \\ 
& $=$ & $\chi (a\otimes (h\otimes f))(k\otimes \widetilde{a}),$%
\end{tabular}
\end{equation*}
i.e. $\chi =\epsilon \circ \alpha .$ Consequently $\chi $ is injective if
and only if $\alpha $ is so.$\blacksquare $
\end{Beweis}

\begin{punto}
\label{comp}Let $H$ be a Hopf $R$-algebra with twisted antipode, $%
A\#_{\sigma }H$ a right $H$-crossed product with invertible cocycle and
consider the $R$-linear maps $\varphi ,\psi :H\otimes _{R}A\rightarrow 
\mathrm{Hom}_{R}(H,A)$ defined as: 
\begin{equation*}
\begin{tabular}{llll}
$\varphi (h\otimes a)(\widetilde{h})$ & $=$ & $\sum [\overline{S}(\widetilde{%
h}_{2})a]\sigma (\overline{S}(\widetilde{h}_{1})\otimes h),$ &  \\ 
$\psi (h\otimes a)(\widetilde{h})$ & $=$ & $\sum \sigma ^{-1}(\widetilde{h}%
_{3}\otimes \overline{S}(\widetilde{h}_{2}))[\widetilde{h}_{4}a]\sigma (%
\widetilde{h}_{5}\otimes \overline{S}(\widetilde{h}_{1})h).$ & 
\end{tabular}
\end{equation*}
Let $U\subseteq H^{\ast }$ be a \emph{right} $H$-module subalgebra, $%
V\subseteq H^{\ast }$ an $R$-submodule and consider the canonical $R$-linear
map $J:A\otimes _{R}V\rightarrow \mathrm{Hom}_{R}(H,A).$ We say $(V,U)$ is 
\emph{compatible}, if the following conditions are satisfied:

\begin{enumerate}
\item  $\varphi (H\otimes _{R}A),$ $\psi (H\otimes _{R}A)\subseteq
J(A\otimes _{R}V);$

\item  $(V,U)$ satisfies the \emph{RL-condition} (\ref{RL}) with respect to $%
H.$
\end{enumerate}
\end{punto}

In the light of Lemma \ref{alp=bet} and the modified RL-condition (\ref{RL})
we introduce an improved version of \cite[Theorem 3, Corollary 4]{Che93}
over arbitrary commutative ground rings:

\begin{proposition}
\label{Chen}Let $H$ be a Hopf $R$-algebra with twisted antipode $\overline{S}%
,$ $A\#_{\sigma }H$ a right $H$-crossed product with invertible cocycle, $%
U\subseteq H^{\ast }$ a right $H$-module subalgebra and consider the $R$%
-pairing $P:=(H,U).$ Assume there exists a \emph{right }$H$\emph{-submodule} 
$V\subseteq H^{\ast },$ such that $(V,U)$ is compatible. If the canonical $R$%
-linear map $\alpha :=\alpha _{A\otimes _{R}H}^{P}:(A\otimes _{R}H)\otimes
_{R}U\rightarrow \mathrm{Hom}_{R}(H,A\otimes _{R}H)$ is injective, then
there exists an $R$-algebra isomorphism 
\begin{equation*}
(A\#_{\sigma }H)\#U\simeq A\otimes _{R}(H\#U).
\end{equation*}
\end{proposition}

\begin{Beweis}
Replacing the inverse of the antipode in \cite[Lemma 2]{Che93} with the
twisted antipode $\overline{S}$ we have a commutative diagram of $R$-algebra
morphisms 
\begin{equation}
\xymatrix{ & (A \#_{\sigma} H) \# U \ar[dl]_{\alpha} \ar[dr]^{\gamma } & \\
\# (H,A \#_{\sigma} H) \ar[rr]^{\pi} & & {\rm End}_{-A} (H \otimes_{R} A) \\
& A \otimes_{R} (H \# U) \ar[ur]_{\chi} \ar[ul]^{\delta } & }  \label{alg-d}
\end{equation}
where 
\begin{equation*}
\begin{tabular}{lll}
$\alpha (a\#(h\#f))(k)$ & $=$ & $(a\#h)f(k),$ \\ 
$\chi (a\otimes (h\#f))(k\otimes \widetilde{a})$ & $=$ & $h(f\rightharpoonup
k)\otimes a\widetilde{a},$ \\ 
$\gamma ((a\#h)\#f)(k\otimes \widetilde{a})$ & $=$ & $\sum
h_{4}(f\rightharpoonup k_{3})\otimes \lbrack \overline{S}(h_{3}k_{2})a]%
\sigma (\overline{S}(h_{2}k_{1})\otimes h_{1})\widetilde{a},$ \\ 
$\delta (a\otimes (h\#f))(k)$ & $=$ & $\sum \sigma ^{-1}(h_{2}k_{4}\otimes 
\overline{S}(h_{1}k_{3}))[(h_{3}k_{5})a]\sigma (h_{4}k_{6}\otimes \overline{S%
}(k_{2}))\#h_{5}(f\rightharpoonup k_{7})\overline{S}(k_{1}),$ \\ 
$\pi (g)(k\otimes \widetilde{a})$ & $=$ & $\nu (\sum g(k_{5})(\sigma
^{-1}(k_{2}\otimes \overline{S}(k_{1}))(k_{3}\rightharpoonup \widetilde{a}%
)\#k_{4}),$%
\end{tabular}
\end{equation*}
and 
\begin{equation*}
\nu :A\#_{\sigma }H\rightarrow H\otimes _{R}A,\text{ }a\#_{\sigma }h\mapsto
\sum h_{4}\otimes \lbrack \overline{S}(h_{3})a]\sigma (\overline{S}%
(h_{2})\otimes h_{1}).
\end{equation*}
Analog to the proof of \cite[Corollary 4]{Che93} (and replacing the inverse
of the antipode by the twisted antipode $\overline{S}$) shows that the
compatibility of $(V,U)$ implies $\mathrm{\func{Im}}(\gamma )\subseteq 
\mathrm{\func{Im}}(\chi )$ and $\mathrm{\func{Im}}(\delta )\subseteq \mathrm{%
\func{Im}}(\alpha ).$ Assume now that $\alpha :=\alpha _{A\otimes _{R}H}^{P}$
is injective. Then $\chi $ is injective by Lemma \ref{alp=bet} and
consequently $\delta $ is injective. Analog to \cite[Lemma 1, Page 2890]
{Che93} $\pi $ is an $R$-algebra isomorphism, hence $\gamma $ is injective
and we are done.$\blacksquare $
\end{Beweis}

\begin{lemma}
\label{M-pure}Let $R$ be Noetherian, $W$ an $R$-module, $U\subseteq W^{\ast }
$ an $R$-submodule and consider the $R$-pairing $P:=(W,U).$ Then the
canonical map $\alpha _{M}^{P}:M\otimes _{R}U\rightarrow \mathrm{Hom}%
_{R}(W,M)$ is injective for an $R$-module $M$ if and only if $U\subset R^{W}$
is $M$-pure. Consequently $P$ satisfies the $\alpha $-condition if and only
if $U\subset R^{W}$ is pure.
\end{lemma}

\begin{Beweis}
Let $M$ be an $R$-module and consider the commutative diagram 
\begin{equation*}
\xymatrix{M \otimes_R U \ar[d]_{id_M \otimes \iota_U} \ar[rr]^{\alpha _M ^P}
& & {\rm Hom}_R (W,M) \ar@{^{(}->}[d] \\ M \otimes_R R^W \ar[rr]_{\varpi} &
& M^W}
\end{equation*}
where $\varpi (m\otimes f)(w)=mf(w).$ Write $M={\underrightarrow{lim}}%
_{I}M_{i}$ as a direct limit of its finitely generated $R$-submodules. Since 
$M_{i}$ is f.p. in $\mathcal{M}_{R}$ we have for every $i\in I$ the
isomorphism of $R$-modules 
\begin{equation*}
\varpi _{i}:M_{i}\otimes R^{W}\rightarrow M_{i}^{W},\text{ }m\otimes
f\mapsto \lbrack w\mapsto mf(w)].
\end{equation*}
Moreover for every $i\in I$ the restriction of $\varpi $ on $M_{i}$
coincides with $\varpi _{i},$ hence 
\begin{equation*}
\varpi =\underrightarrow{lim}\varpi _{M_{i}}:\underrightarrow{lim}%
M_{i}\otimes R^{W}\rightarrow \underrightarrow{lim}M_{i}^{W}\subset M^{W}
\end{equation*}
is injective. It's obvious then that $\alpha _{M}^{P}$ is injective iff $%
id_{M}\otimes \iota _{U}$ is injective iff $U\subset R^{W}$ is $M$-pure.$%
\blacksquare $
\end{Beweis}

We are ready now to present the main duality theorem in this note, which
generalizes \cite[Corollary 4]{Che93} (resp. \cite[Corollary 5]{Che93}) from
the case of a base field (resp. a Dedekind domain) to the case of an
arbitrary Noetherian ring:

\begin{theorem}
\label{haupt}Let $R$ be Noetherian, $H$ a Hopf $R$-algebra with twisted
antipode, $A\#_{\sigma }H$ a right $H$-crossed product with invertible
cocycle, $U\subseteq H^{\ast }$ a right $H$-module subalgebra and consider
the $R$-pairing $P:=(H,U).$ Assume there exists a right $H$-submodule $%
V\subseteq H^{\ast },$ such that $(V,U)$ is compatible. If $U\subset R^{H}$
is $A\otimes _{R}H$-pure \emph{(}e.g. $H$ is a Hopf $\alpha $-algebra and $%
U\subseteq H^{\circ }$ is an $R$-subbialgebra\emph{)}, then we have an $R$%
-algebra isomorphism 
\begin{equation*}
(A\#_{\sigma }H)\#U\simeq A\otimes _{R}(H\#U).
\end{equation*}
\end{theorem}

\begin{Beweis}
By Proposition \ref{Chen} it remains to show that $\alpha _{A\otimes
_{R}H}^{P}$ is injective. If $U\subset R^{H}$ is $A\otimes _{R}H$-pure, then 
$\alpha _{A\otimes _{R}H}^{P}$ is injective by Lemma \ref{M-pure}. If $H$ is
a Hopf $\alpha $-algebra, then $H^{\circ }\subset R^{H}$ is pure and for
every $R$-subbialgebra $U\subseteq H^{\circ },$ $U\subset R^{H}$ is pure
(since by convention $U\subseteq H^{\circ }$ is pure), hence $\alpha
_{A\otimes _{R}H}^{P}$ is injective.$\blacksquare $
\end{Beweis}

\begin{definition}
Let $R$ be Noetherian. After \cite{Mon93} we call an $R$-algebra $A$ \emph{%
residually finite }(called in other references \emph{proper})\emph{, }if the
canonical map $A\rightarrow A^{\circ \ast }$ is injective (equivalently, if $%
\bigcap \{\mathrm{Ke}(f)\mid f\in A^{\circ }\}=0$).
\end{definition}

\begin{corollary}
\label{dense}Let $H$ be a Hopf $R$-algebra with twisted antipode and $_{R}H$
projective, $A\#_{\sigma }H$ a right $H$-crossed product with invertible
cocycle, $U\subseteq H^{\ast }$ a right $H$-module subalgebra and consider
the $R$-paring $P:=(H,U).$ Assume there exists a right $H$-submodule $%
V\subseteq H^{\ast },$ such that $(V,U)$ is compatible. If $U\subseteq
H^{\ast }$ is dense and the canonical $R$-linear map $\alpha _{A\otimes
_{R}H}^{P}$ is injective \emph{(}e.g. $R$ is Noetherian and $U\subseteq R^{H}
$ is $A$-pure\emph{)}, then there exists a dense $R$-subalgebra $\mathcal{L}%
\subseteq \mathrm{End}_{R}(H)$ and an $R$-algebra isomorphism 
\begin{equation*}
(A\#_{\sigma }H)\#U\simeq A\otimes _{R}\mathcal{L}.
\end{equation*}
This is the case in particular, if $R$ is a QF ring, $H$ is a residually
finite Hopf $\alpha $-algebra and $U\subseteq H^{\circ }$ is a dense $R$%
-subbialgebra.
\end{corollary}

\begin{Beweis}
If $U\subseteq H^{\ast }$ is dense, then $\mathcal{L}:=H\#U\overset{\lambda 
}{\hookrightarrow }\mathrm{End}_{R}(H)$ is a dense $R$-subalgebra by Lemma 
\ref{density} (4) and the isomorphism follows by Theorem \ref{haupt}. If $R$
is a QF ring and $H$ is a residually finite Hopf $\alpha $-algebra, then $%
H^{\circ }\subset H^{\ast }$ is dense by \cite[Proposition 2.4.19]{Abu2001}.
If moreover $U\subseteq H^{\circ }$ is a dense $R$-subbialgebra, then $%
U\subseteq H^{\ast }$ is dense, $\alpha _{A\otimes _{R}H}^{P}$ is injective
and we are done.$\blacksquare $
\end{Beweis}

\begin{corollary}
\label{f.d.}Let $H$ be a Hopf $R$-algebra with twisted antipode, $%
A\#_{\sigma }H$ a right $H$-crossed product with invertible cocycle and
consider the $R$-pairing $P:=(H,H^{\ast }).$ Then we have an isomorphism of $%
R$-algebras 
\begin{equation*}
(A\#_{\sigma }H)\#H^{\ast }\simeq A\otimes _{R}(H\#H^{\ast }).
\end{equation*}
at least when:

\begin{enumerate}
\item  $_{R}H$ is finitely generated projective, \emph{or}

\item  $_{R}A$ is finitely generated, $H$ is cocommutative and $\alpha
_{A\otimes _{R}H}^{P}$ is injective \emph{(}e.g. $R$ is Noetherian and $%
H^{\ast }\hookrightarrow R^{H}$ is $A\otimes _{R}H$-pure\emph{)}.
\end{enumerate}
\end{corollary}

\begin{Beweis}
\begin{enumerate}
\item  Since $_{R}H$ is finitely generated projective, the canonical $R$%
-linear map $J:A\otimes _{R}H^{\ast }\rightarrow \mathrm{Hom}_{R}(H,A)$ is
bijective and $H^{\ast }$ satisfies the RL-condition (\ref{RL}) with respect
to $H$ by Lemma \ref{density} (3), hence $(H^{\ast },H^{\ast })$ is
compatible. Moreover $P=(H,H^{\ast })\simeq (H^{\ast \ast },H^{\ast })$
satisfies the $\alpha $-condition, since $_{R}H^{\ast }$ is finitely
generated projective. The result follows now by Proposition \ref{Chen}.

\item  Since $_{R}A$ is finitely generated, the canonical $R$-linear map $%
J:A\otimes _{R}H^{\ast }\rightarrow \mathrm{Hom}_{R}(H,A)$ is surjective.
Since $H$ is cocommutative, $H^{\ast }$ satisfies the RL-condition (\ref{RL}%
) with respect to $H,$ hence $(H^{\ast },H^{\ast })$ is compatible. By
assumption $\alpha _{A\otimes _{R}H}^{P}$ is injective and we are done by
Proposition \ref{Chen}.$\blacksquare $
\end{enumerate}
\end{Beweis}

\begin{corollary}
\label{Mn}Let $H$ be a free Hopf $R$-algebra of rank $n$ and $A\#_{\sigma }H$
a right $H$-crossed product with invertible cocycle. Then we have an
isomorphism of $R$-algebras 
\begin{equation*}
(A\#_{\sigma }H)\#H^{\ast }\simeq A\otimes _{R}M_{n}(R)\simeq M_{n}(A).
\end{equation*}
\end{corollary}

\begin{Beweis}
By Corollary \ref{f.d.} $(A\#_{\sigma }H)\#H^{\ast }\simeq A\otimes
_{R}(H\#H^{\ast }).$ Since $_{R}H$ is finitely generated projective, $%
H\#H^{\ast }\simeq \mathrm{End}_{R}(H)$ by Lemma \ref{density} (3). But $%
_{R}H$ is free of rank $n,$ hence $\mathrm{End}_{R}(H)\simeq M_{n}(R).$ It's
evident that $A\otimes _{R}M_{n}(R)\simeq M_{n}(A)$ and we are done.$%
\blacksquare $
\end{Beweis}

\section*{The right $H^{\ast }$-submodule $H^{\protect\upsilon }\subseteq
H^{\ast }$}

In what follows let $H$ be a \emph{locally projective} Hopf $R$-algebra with 
\emph{twisted antipode} and consider the measuring $\alpha $-pairing $%
P:=(H^{\ast },H)$ (notice that the canonical $R$-linear map $\alpha
_{R}^{P}:H\rightarrow H^{\ast \ast }$ is injective).

\begin{lemma}
Consider $H^{\ast }$ with the right $H^{\ast }$-action 
\begin{equation*}
(f\leftharpoonup g)(h):=\sum g(h_{3}\overline{S}(h_{1}))f(h_{2})\text{ for
all }f,g\in H^{\ast }\text{ and }h\in H.
\end{equation*}
Then $H^{\ast }$ is a right $H^{\ast }$-module and $H^{\upsilon }:=$ $^{H}%
\mathrm{Rat}(H_{H^{\ast }}^{\ast })$ is a left $H$-comodule with structure
map $\upsilon :H^{\upsilon }\rightarrow H\otimes _{R}H^{\upsilon }.$
\end{lemma}

\begin{Beweis}
For arbitrary $f,g,\widetilde{g}\in H^{\ast }$ we have 
\begin{equation*}
\begin{tabular}{lll}
$(f\leftharpoonup (g\star \widetilde{g}))(h)$ & $=$ & $\sum (g\star 
\widetilde{g})(h_{3}\overline{S}(h_{1}))f(h_{2})$ \\ 
& $=$ & $\sum g(h_{31}\overline{S}(h_{1})_{1})\widetilde{g}(h_{32}\overline{S%
}(h_{1})_{2})f(h_{2})$ \\ 
& $=$ & $\sum g(h_{31}\overline{S}(h_{12}))\widetilde{g}(h_{32}\overline{S}%
(h_{11}))f(h_{2})$ \\ 
& $=$ & $\sum g(h_{4}\overline{S}(h_{2}))\widetilde{g}(h_{5}\overline{S}%
(h_{1}))f(h_{3})$ \\ 
& $=$ & $\sum \widetilde{g}(h_{3}\overline{S}(h_{1}))g(h_{23}\overline{S}%
(h_{21}))f(h_{22})$ \\ 
& $=$ & $\sum \widetilde{g}(h_{3}\overline{S}(h_{1}))(f\leftharpoonup
g)(h_{2})$ \\ 
& $=$ & $((f\leftharpoonup g)\leftharpoonup \widetilde{g})(h).$%
\end{tabular}
\end{equation*}
Since $_{R}H$ is locally projective, we have analog to Theorem \ref{equal}
that $^{H}\mathrm{Rat}(H_{H^{\ast }}^{\ast })$ is a left $H$-comodule.$%
\blacksquare $
\end{Beweis}

\begin{proposition}
\label{Hw-prop}Consider the left $H$-comodule $(H^{\upsilon },\upsilon ).$

\begin{enumerate}
\item  If $f\in H^{\upsilon },$ then $\upsilon (f)=\sum f_{<-1>}\otimes
f_{<0>}$ satisfies the following conditions:

\begin{enumerate}
\item  $f\star g=\sum gf_{<-1>}\star f_{<0>}$ for all $g\in H^{\ast }.$

\item  $h\leftharpoonup f=\sum f_{<-1>}(f_{<0>}\rightharpoonup h)$ for all $%
h\in H.$

\item  $\sum h_{3}\overline{S}(h_{1})f(h_{2})=\sum f_{<-1>}f_{<0>}(h)$ for
all $h\in H.$
\end{enumerate}

\item  Let $f\in H^{\ast }.$ If there exists $\zeta =\sum f_{<-1>}\otimes
f_{<0>}\in H\otimes _{R}H^{\ast }$ that satisfies any of the conditions in 
\emph{(1)}, then $f\in H^{\upsilon }$ and $\upsilon (f)=\zeta .$

\item  For all $f,\widetilde{f}\in H^{\upsilon }$ and $g\in H^{\ast }$ we
have 
\begin{equation*}
(f\star \widetilde{f})\star g=\sum g(\widetilde{f}_{<-1>}f_{<-1>})\star
(f_{<0>}\star \widetilde{f}_{<0>}).
\end{equation*}

\item  $H^{\upsilon }\subseteq H^{\ast }$ is a right $H$-submodule with 
\begin{equation*}
\upsilon (fh)=\sum \overline{S}(h_{3})f_{<-1>}h_{1}\otimes f_{<0>}h_{2}\text{
for all }h\in H\text{ and }f\in H^{\upsilon }.
\end{equation*}
\end{enumerate}
\end{proposition}

\begin{Beweis}
\begin{enumerate}
\item  Let $f\in H^{\upsilon }$ with $\upsilon (f)=\sum f_{<-1>}\otimes
f_{<0>}.$

\begin{enumerate}
\item  For all $g\in H^{\ast }$ and $h\in H$ we have 
\begin{equation*}
\begin{tabular}{lllll}
$(f\star g)(h)$ & $=$ & $\sum f(h_{1})g(h_{2})$ & $=$ & $\sum g(h_{3}%
\overline{S}(h_{12})h_{11})f(h_{2})$ \\ 
& $=$ & $\sum g(h_{23}\overline{S}(h_{21})h_{1})f(h_{22})$ & $=$ & $\sum
(h_{1}g)(h_{23}\overline{S}(h_{21}))f(h_{22})$ \\ 
& $=$ & $\sum (f\leftharpoonup (h_{1}g))(h_{2})$ & $=$ & $\sum
(h_{1}g)(f_{<-1>})f_{<0>}(h_{2})$ \\ 
& $=$ & $\sum g(f_{<-1>}h_{1})f_{<0>}(h_{2})$ & $=$ & $\sum
(gf_{<-1>})(h_{1})f_{<0>}(h_{2})$ \\ 
& $=$ & $(\sum (gf_{<-1>})\star f_{<0>})(h).$ &  & 
\end{tabular}
\end{equation*}

\item  For all $g\in H^{\ast }$ and $h\in H$ we have 
\begin{equation*}
\begin{tabular}{lllll}
$g(h\leftharpoonup f)$ & $=$ & $g(\sum f(h_{1})h_{2})$ & $=$ & $\sum
f(h_{1})g(h_{2})$ \\ 
& $=$ & $\sum f(h_{2})g(h_{3}\overline{S}(h_{12})h_{11})$ & $=$ & $\sum
f(h_{22})g(h_{23}\overline{S}(h_{21})h_{1})$ \\ 
& $=$ & $\sum f(h_{22})(h_{1}g)(h_{23}\overline{S}(h_{21}))$ & $=$ & $\sum
(f\leftharpoonup (h_{1}g))(h_{2})$ \\ 
& $=$ & $\sum (h_{1}g)(f_{<-1>})f_{<0>}(h_{2})$ & $=$ & $g(\sum
f_{<-1>}h_{1}f_{<0>}(h_{2}))$ \\ 
& $=$ & $g(\sum f_{<-1>}(f_{<0>}\rightharpoonup h)).$ &  & 
\end{tabular}
\end{equation*}

\item  Trivial.
\end{enumerate}

\item  Let $f\in H^{\ast }$ and $\zeta =\sum f_{<-1>}\otimes f_{<0>}\in
H\otimes _{R}H^{\ast }.$ We are done once we have shown that $%
(f\leftharpoonup g)(h)=\sum g(f_{<-1>})f_{<0>}(h)$ for arbitrary $g\in
H^{\ast }$ and $h\in H.$

\begin{enumerate}
\item  Assume (1-a) holds. Then we have 
\begin{equation*}
\begin{tabular}{lll}
$(f\leftharpoonup g)(h)$ & $=$ & $\sum g(h_{3}\overline{S}(h_{1}))f(h_{2})$
\\ 
& $=$ & $\sum (\overline{S}(h_{1})g)(h_{22})f(h_{21})$ \\ 
& $=$ & $\sum (f\star \overline{S}(h_{1})g)(h_{2})$ \\ 
& $=$ & $\sum (\overline{S}(h_{1})gf_{<-1>}\star f_{<0>})(h_{2})$ \\ 
& $=$ & $\sum ((\overline{S}(h_{1})gf_{<-1>})(h_{2})f_{<0>}(h_{3})$ \\ 
& $=$ & $\sum g(f_{<-1>}h_{2}\overline{S}(h_{1}))f_{<0>}(h_{3})$ \\ 
& $=$ & $\sum g(f_{<-1>})f_{<0>}(h).$%
\end{tabular}
\end{equation*}

\item  Assume (1-b) holds. Then we have 
\begin{equation*}
\begin{tabular}{lll}
$(f\leftharpoonup g)(h)$ & $=$ & $\sum g(h_{3}\overline{S}(h_{1}))f(h_{2})$
\\ 
& $=$ & $\sum (\overline{S}(h_{1})g)(h_{2}\leftharpoonup f)$ \\ 
& $=$ & $\sum (\overline{S}(h_{1})g)(f_{<-1>}(f_{<0>}\rightharpoonup h_{2}))$
\\ 
& $=$ & $\sum (\overline{S}(h_{1})g)(f_{<-1>}h_{2}f_{<0>}(h_{3}))$ \\ 
& $=$ & $\sum g(f_{<-1>}h_{2}\overline{S}(h_{1}))f_{<0>}(h_{3})$ \\ 
& $=$ & $\sum g(f_{<-1>})f_{<0>}(h).$%
\end{tabular}
\end{equation*}

\item  Trivial.
\end{enumerate}

\item  Let $f,\widetilde{f}\in H^{\upsilon }.$ For arbitrary $g\in H^{\ast }$
we have by (1-a): 
\begin{equation*}
\begin{tabular}{lll}
$(f\star \widetilde{f})\star g$ & $=$ & $f\star (\widetilde{f}\star g)$ \\ 
& $=$ & $\sum f\star (g\widetilde{f}_{<-1>}\star \widetilde{f}_{<0>})$ \\ 
& $=$ & $\sum (f\star g\widetilde{f}_{<-1>})\star \widetilde{f}_{<0>}$ \\ 
& $=$ & $\sum (g\widetilde{f}_{<-1>})f_{<-1>}\star (f_{<0>}\star \widetilde{f%
}_{<0>})$ \\ 
& $=$ & $\sum g(\widetilde{f}_{<-1>}f_{<-1>})\star (f_{<0>}\star \widetilde{f%
}_{<0>})$%
\end{tabular}
\end{equation*}

\item  Let $f\in H^{\upsilon }$ and $h\in H.$ Then we have for all $g\in
H^{\ast }$ and $k\in H:$%
\begin{equation*}
\begin{tabular}{lll}
$((fh)\leftharpoonup g)(k)$ & $=$ & $\sum g(k_{3}\overline{S}%
(k_{1}))(fh)(k_{2})$ \\ 
& $=$ & $\sum g(k_{3}\overline{S}(k_{1}))f(hk_{2})$ \\ 
& $=$ & $\sum g(\overline{S}(h_{4})h_{3}k_{3}\overline{S}(k_{1})\overline{S}%
(h_{12})h_{11})f(h_{2}k_{2})$ \\ 
& $=$ & $\sum g(\overline{S}(h_{3})h_{23}k_{3}\overline{S}(k_{1})\overline{S}%
(h_{21})h_{1})f(h_{22}k_{2})$ \\ 
& $=$ & $\sum (h_{1}g\overline{S}(h_{3}))(h_{23}k_{3}\overline{S}(k_{1})%
\overline{S}(h_{21}))f(h_{22}k_{2})$ \\ 
& $=$ & $\sum (h_{1}g\overline{S}(h_{3}))(h_{23}k_{3}\overline{S}%
(h_{21}k_{1}))f(h_{22}k_{2})$ \\ 
& $=$ & $\sum (h_{1}g\overline{S}(h_{3}))((h_{2}k)_{3}\overline{S}%
((h_{2}k)_{1}))f((h_{2}k)_{2})$ \\ 
& $=$ & $\sum (f\leftharpoonup (h_{1}g\overline{S}(h_{3})))(h_{2}k)$ \\ 
& $=$ & $\sum (h_{1}g\overline{S}(h_{3}))(f_{<-1>})f_{<0>}(h_{2}k)$ \\ 
& $=$ & $\sum g(\overline{S}(h_{3})f_{<-1>}h_{1})(f_{<0>}h_{2})(k).%
\blacksquare $%
\end{tabular}
\end{equation*}
\end{enumerate}
\end{Beweis}

\qquad As a consequence of Proposition \ref{Chen} and Theorem \ref{haupt}
and analog to \cite[Theorem 4.2]{Kop92} we get

\begin{theorem}
\label{v-hm1}Let $H$ be a locally projective Hopf $R$-algebra with twisted
antipode, $U\subseteq H^{\ast }$ a right $H$-module subalgebra, $P:=(H,U)$
the induced $R$-pairing and assume that $\varphi (H\otimes _{R}A),$ $\psi
(H\otimes _{R}A)\subseteq J(A\otimes _{R}\upsilon ^{-1}(H\otimes _{R}U)).$
If $\alpha _{A\otimes _{R}H}^{P}$ is injective (e.g. $R$ is Noetherian and $%
U\subset R^{H}$ is $A$-pure), then there is an $R$-algebra isomorphism 
\begin{equation*}
(A\#_{\sigma }H)\#U\simeq A\otimes _{R}(H\#U).
\end{equation*}
\end{theorem}

\begin{Beweis}
It follows from Proposition \ref{Hw-prop} (4) that $V:=\upsilon
^{-1}(H\otimes _{R}U)\subset H^{\ast }$ is a right $H$-submodule. Since $%
V\subseteq H^{\upsilon },$ it follows by Proposition \ref{Hw-prop} (1-b)
that $(V,U)$ satisfies the RL-condition (\ref{RL}) with respect to $H.$
Consequently $(V,U)$ is compatible. If $\alpha _{A\otimes _{R}H}^{P}$ is
injective, then the result follows by Proposition \ref{Chen}.$\blacksquare $
\end{Beweis}

\begin{corollary}
\label{norm}Let $H$ be a locally projective Hopf $R$-algebra with twisted
antipode, $U\subseteq H^{\upsilon }$ a right $H$-module subalgebra of $%
H^{\ast },$ $P:=(H,U)$ the induced $R$-pairing and assume that $\upsilon
(U)\subseteq H\otimes _{R}U$ and $\overline{\varphi }(H\otimes _{R}A),$ $%
\overline{\psi }(H\otimes _{R}A)\subseteq J(A\otimes _{R}U).$ If $\alpha
_{A\otimes _{R}H}^{P}$ is injective (e.g. $R$ is Noetherian and $U\subset
R^{H}$ is $A$-pure), then there is an $R$-algebra isomorphism 
\begin{equation*}
(A\#_{\sigma }H)\#U\simeq A\otimes _{R}(H\#U).
\end{equation*}
\end{corollary}

\section*{Blattner-Montgomery Duality Theorem revisited}

The following definition is suggested by \cite[Definition 1.3]{BM85}:

\begin{definition}
Let $R$ be Noetherian, $H$ an $R$-bialgebra, $U\subseteq H^{\circ }$ an $R$%
-submodule and $A$ a left $H$-module algebra. Then $A$ will be called $U$%
\textbf{-locally finite }if and only if for every $a\in A$ there exists a
finite subset $\{f_{1},...,f_{k}\}\subset U,$ such that $\bigcap_{i=1}^{k}%
\mathrm{Ke}(f_{i})\subseteq (0_{A}:a).$
\end{definition}

\begin{lemma}
\emph{(\cite[Proposition 3.3]{Abu-a})}\label{cc-ma} Let $R$ be Noetherian, $H
$ an $\alpha $-bialgebra, $U\subseteq H^{\circ }$ an $R$-subbialgebra and
consider the measuring $\alpha $-pairing $(H,U).$

\begin{enumerate}
\item  If $A$ is a right \emph{(}a\emph{\ }left\emph{)} $U$-comodule
algebra, then $A$ is a left \emph{(}a right\emph{) }$H$-module algebra.

\item  If $A$ is a left \emph{(}a right\emph{) }$H$-module algebra, then $%
\mathrm{Rat}^{U}(_{H}A)$ is a right \emph{(}a left\emph{)} $U$-comodule
algebra.
\end{enumerate}
\end{lemma}

\qquad The following result generalizes \cite[Lemma 1.5]{BM85} from the case
of a base field to the case of an arbitrary Noetherian ground ring.

\begin{lemma}
\label{U-loc}Let $R$ be Noetherian, $A$ an $R$-algebra, $H$ an $\alpha $%
-bialgebra and $U\subseteq H^{\circ }$ an $R$-subbialgebra. Then $A$ is a $U$%
-locally finite left $H$-module algebra if and only if $A$ is a right $U$%
-comodule algebra.
\end{lemma}

\begin{Beweis}
Consider the measuring $\alpha $-pairing $(H,U).$ Assume $A$ to be a right $U
$-comodule algebra. Then $A$ is a left $H$-module algebra by Lemma \ref
{cc-ma} (1). Moreover for every $a\in A$ with $\varrho
(a)=\sum\limits_{j=1}^{n}a_{j}\otimes g_{j}\in A\otimes _{R}U$ we have $%
\bigcap_{j=1}^{n}\mathrm{Ke}(g_{i})\subseteq (0_{A}:a),$ i.e. $_{H}A$ is $U$%
-locally finite. On the other hand, assume $A$ to be a $U$-locally finite
left $H$-module algebra and consider $H$ with the left $U$-adic topology $%
\mathcal{T}_{U-}(H)$ (see \ref{adic-top}). By Lemma \ref{cc-ma} (2) $\mathrm{%
Rat}^{U}(_{H}A)$ is a right $U$-comodule algebra and we are done once we
have shown that $A=\mathrm{Rat}^{U}(_{H}A).$ By assumption there exists for
every $a\in A$ a subset $W=\{f_{1},...,f_{k}\}\subset U,$ such that $%
\bigcap_{i=1}^{k}\mathrm{Ke}(f_{i})\subseteq (0_{A}:a).$ If $h\in (0_{U}:W),$
then $f_{i}(h)=(hf_{i})(1_{H})=0$ for $i=1,...,k$ and so $(0_{U}:W)\subseteq
\bigcap_{i=1}^{k}\mathrm{Ke}(f_{i})\subseteq (0_{A}:a),$ i.e. $A$ is a
discrete left $(H,\mathcal{T}_{U-}(H))$-module. Consequently $A\in \sigma
\lbrack _{H}U]=\mathrm{Rat}^{U}(_{H}\mathcal{M})$ (see Theorem \ref{equal}),
i.e. $A=\mathrm{Rat}^{U}(_{H}A).\blacksquare $
\end{Beweis}

The following result provides an improved version of Blattner-Montgomery
duality theorem for the case of arbitrary Noetherian base rings, replacing
the assumption ``$U\subseteq H^{\circ }$ is a Hopf $R$-subalgebra with
bijective antipode'' in the original version \cite[Theorem 2.1]{BM85} and in 
\cite[3.2]{AG-TL2001} with ``$U\subseteq H^{\circ }$ is any $R$%
-subbialgebra'' (as suggested by M. Koppinen \cite[Corollary 5.4]{Kop92});
and replacing the assumption that $H$ has a bijective antipode with the
weaker condition that $H$ has a twisted antipode $\overline{S}.$

\begin{corollary}
\label{BM}Let $R$ be Noetherian, $H$ a Hopf $\alpha $-algebra with twisted
antipode and $U\subseteq $ $H^{\circ }$ an $R$-subbialgebra.{\normalsize \ }%
Let $A$ be a $U$-locally finite left $H$-module algebra and consider $A$
with the induced right $H$-comodule structure. If there exists a right $H$%
-submodule $V\subseteq H^{\ast },$ such that $\mathrm{Cf}(A)\cup \overline{S}%
^{\ast }(\mathrm{Cf}(A))\subseteq V$ and $(V,U)$ satisfies the RL-condition 
\emph{(\ref{RL}) }with respect to $H,$ then we have an isomorphism of $R$%
-algebras 
\begin{equation*}
(A\#H)\#U\simeq A\otimes _{R}(H\#U).
\end{equation*}
\end{corollary}

\begin{Beweis}
For the trivial cocycle $\sigma (h\otimes k):=\varepsilon (h)\varepsilon
(k)1_{A}$ we have $A\#_{\sigma }H=A\#H.$ Consider the canonical $R$-linear
map $J:A\otimes _{R}V\rightarrow \mathrm{Hom}_{R}(H,A).$ For every $h\in H$
and $a\in A$ we have 
\begin{equation*}
\begin{tabular}{lllll}
$\varphi (h\otimes a)(\widetilde{h})$ & $=$ & $[\overline{S}(\widetilde{h}%
)a]\varepsilon (h)$ & $=$ & $\sum a_{<0>}<\overline{S}(\widetilde{h}%
),a_{<1>}>\varepsilon (h)$ \\ 
& $=$ & $\sum a_{<0>}\overline{S}^{\ast }(a_{<1>})(\widetilde{h})\varepsilon
(h)$ & $=$ & $J(\sum a_{<0>}\varepsilon (h)\otimes \overline{S}^{\ast
}(a_{<1>}))(\widetilde{h}),$%
\end{tabular}
\end{equation*}
i.e. $\varphi (H\otimes _{R}A)\subseteq J(A\otimes _{R}V).$ On the other
hand we have for all $h,$ $\widetilde{h}\in H$ and $a\in A:$%
\begin{equation*}
\psi (h\otimes a)(\widetilde{h})=[\widetilde{h}a]\varepsilon (h)=\sum
a_{<0>}<\widetilde{h},a_{<1>}>\varepsilon (h)=J(\sum a_{<0>}\varepsilon
(h)\otimes a_{<1>})(\widetilde{h}),
\end{equation*}
i.e. $\psi (H\otimes _{R}A)\subseteq J(A\otimes _{R}U).$ By assumption $(V,U)
$ satisfies the RL-condition (\ref{RL}) with respect to $H,$ hence $(V,U)$
is compatible and the result follows then by Theorem \ref{haupt} (notice
that $P=(H,U)$ is an $\alpha $-pairing).$\blacksquare $
\end{Beweis}

As a consequence of Corollaries \ref{dense} and \ref{BM} we get

\begin{corollary}
\label{BM-dense}Let $R$ be Noetherian, $H$ a projective Hopf $\alpha $%
-algebra with twisted antipode and $U\subseteq $ $H^{\circ }$ an $R$%
-subbialgebra.{\normalsize \ }Let $A$ be a $U$-locally finite left $H$%
-module algebra and consider $A$ with the induced right $H$-comodule
structure. Assume there exists a right $H$-submodule $V\subseteq H^{\ast },$
such that $\mathrm{Cf}(A)\cup \overline{S}^{\ast }(\mathrm{Cf}(A))\subseteq V
$ and $(V,U)$ satisfies the RL-condition \emph{(\ref{RL}) }with respect to $%
H.$ If $U\subseteq H^{\ast }$ is dense, then there exists a dense $R$%
-subalgebra $\mathcal{L}\subseteq \mathrm{End}_{R}(H)\ $and an $R$-algebra
isomorphism 
\begin{equation*}
(A\#H)\#U\simeq A\otimes _{R}\mathcal{L}.
\end{equation*}
In particular this holds, if $R$ is a QF ring, $H$ is residually finite and $%
U\subseteq H^{\circ }$ is dense.
\end{corollary}

\section*{Cleft $H$-extensions}

Hopf-Galois extensions were presented by S. Chase and M. Sweedler \cite{CS69}
for a \emph{commutative }$R$-algebra acting on a Hopf $R$-algebra and are
considered as generalization of the classical Galois extensions over fields
(e.g. \cite[8.1.2]{Mon93}). In \cite{KT81} H. Kreimer and M. Takeuchi
extended these to the \emph{noncommutative }case.

\begin{punto}
$H$\textbf{-Extensions. }(\cite{Doi85}) Let $H$ be an $R$-bialgebra, $B$ a
right $H$-comodule algebra and consider the $R$-algebra $A:=B^{coH}=\{a\in
B\mid \varrho (a)=a\otimes 1_{H}\}.$ The algebra extension $A\hookrightarrow
B$ is called a \emph{right }$H$\emph{-extension. }A (\emph{total}) \emph{%
integral} for $B$ is an $H$\emph{-colinear} map $\theta :H\rightarrow B$
(with $\theta (1_{H})=1_{B}$). If $B$ admits an integral, which is
invertible in $(\mathrm{Hom}_{R}(H,B),\star ),$ then the right $H$-extension 
$A\hookrightarrow B$ is called \emph{cleft.}
\end{punto}

\begin{lemma}
\label{eq}\emph{(\cite[Theorems 9, 11]{DT86}, \cite[Theorem 1.18]{BM89}, 
\cite[1.1.]{DT89}) }Let $H$ be an $R$-bialgebra.

\begin{enumerate}
\item  If $B/A$ is a \emph{cleft }right $H$-extension with total invertible
integral $\theta :H\rightarrow B,$ then $A$ is a left $H$-module algebra
through 
\begin{equation*}
ha=\sum \theta (h_{1})a\theta ^{-1}(h_{2})\text{ for all }h\in H\text{ and }%
a\in A
\end{equation*}
and $A\#_{\sigma }H$ is a right $H$-crossed product with invertible cocycle 
\begin{equation*}
\sigma (h\otimes k)=\sum \theta (h_{1})\theta (k_{1})\theta
^{-1}(h_{2}k_{2}),\text{ where }\sigma ^{-1}(h\otimes k)=\sum \theta
(h_{1}k_{1})\theta ^{-1}(k_{2})\theta ^{-1}(h_{2}).
\end{equation*}
Moreover $B\simeq A\#_{\sigma }H$ as right $H$-comodule algebras.

\item  Let $H$ be a Hopf $R$-algebra. If $B:=A\#_{\sigma }H$ is a right $H$%
-crossed product with invertible cocycle $\sigma \in \mathrm{Hom}%
_{R}(H\otimes _{R}H,A),$ then $B/A$ is a \emph{cleft }right $H$-extension
with invertible total integral 
\begin{equation*}
\theta :H\rightarrow A\#_{\sigma }H,\text{ }\theta (h)=1_{A}\#h,\text{ where 
}\theta ^{-1}(h)=\sum \sigma ^{-1}(S(h_{2})\otimes h_{3})\#_{\sigma
}S(h_{1}).
\end{equation*}
\end{enumerate}
\end{lemma}

\qquad Let $H$ be a Hopf $R$-algebra with twisted antipode, $B/A$ a cleft
right $H$-extension with invertible total integral $\theta :H\rightarrow B$
and consider the $R$-linear maps $\widetilde{\varphi },\widetilde{\psi }%
:A\otimes _{R}H\rightarrow \mathrm{Hom}_{R}(H,A)$ defined as: 
\begin{eqnarray}
\widetilde{\varphi }(h\otimes a)(\widetilde{h}) &=&\sum \theta (\overline{S}(%
\widetilde{h}_{2}))a\theta (h_{1})\theta ^{-1}(\overline{S}(\widetilde{h}%
_{1})h_{2}),  \label{varphi-tel} \\
\widetilde{\psi }(h\otimes a)(\widetilde{h}) &=&\sum \theta ^{-1}(\overline{S%
}(\widetilde{h}_{3}))a\theta (\overline{S}(\widetilde{h}_{2})h_{1})\theta
^{-1}(\widetilde{h}_{4}\overline{S}(\widetilde{h}_{1})h_{2}).
\label{psi-tel}
\end{eqnarray}

With the help of Lemma \ref{eq} one can easily derive the following version
of Theorem \ref{haupt} and Corollary \ref{dense} for \emph{cleft }right $H$%
-extensions:

\begin{theorem}
\label{cleft-du}Let $R$ be Noetherian, $H$ a Hopf $R$-algebra with twisted
antipode, $B/A$ a \emph{cleft} right $H$-extension with invertible total
integral $\theta :H\rightarrow B,$ $U\subseteq H^{\ast }$ a right $H$-module
subalgebra and consider the $R$-pairing $P:=(H,U).$ Assume there exists a
right $H$-submodule $V\subseteq H^{\ast },$ such that:

\begin{enumerate}
\item  $\widetilde{\varphi }(H\otimes _{R}A),$ $\widetilde{\psi }(H\otimes
_{R}A)\subseteq J(A\otimes _{R}V);$

\item  $(V,U)$ satisfies the \emph{RL-condition (\ref{RL})} with respect to $%
H.$
\end{enumerate}

If $U\subset R^{H}$ is $A\otimes _{R}H$-pure \emph{(}e.g. $H$ is a Hopf $%
\alpha $-algebra and $U\subseteq H^{\circ }$ is an $R$-subbialgebra\emph{)},
then there is an $R$-algebra isomorphism 
\begin{equation*}
B\#U\simeq A\otimes _{R}(H\#U).
\end{equation*}
If moreover $_{R}H$ is projective and $U\subseteq H^{\ast }$ is dense \emph{(%
}e.g. $R$ is a QF ring, $H$ is residually finite and $U\subseteq H^{\circ }$
is dense), then $B\#U\simeq A\otimes _{R}\mathcal{L}$ for a \emph{dense} $R$%
-subalgebra $\mathcal{L}\subseteq \mathrm{End}_{R}(H).$
\end{theorem}

\section{Koppinen Duality Theorem}

In this section we prove an improved version of Koppinen's duality theorem
presented in \cite{Kop92} over arbitrary Noetherian ground rings under
fairly weak conditions. In fact the results in this section are similar to
those in the second section with a main advantage, that they are evident for
arbitrary Hopf $R$-algebras (not necessarily with twisted antipodes).

\begin{punto}
\label{smash-op}Let $H$ be an $R$-bialgebra and $B$ a right $H$-comodule
algebra. Then $\#^{op}(H,B)=\mathrm{Hom}_{R}(H,B)$ is an associative $R$%
-algebra with multiplication 
\begin{equation}
(f\widetilde{\star }g)(h)=\sum f(h_{2})_{<0>}g(h_{1}f(h_{2})_{<1>})\text{
for all }f,g\in \mathrm{Hom}_{R}(H,B),\text{ }h\in H  \label{bls}
\end{equation}
and unity $\eta _{B}\circ \varepsilon _{H}.$ If $U\subseteq H^{\ast }$ is a 
\emph{left }$H$-module subalgebra (with $\varepsilon _{H}\in U$), then $%
B\#^{op}U=B\otimes _{R}U$ is an associative $R$-algebra with multiplication 
\begin{equation}
(b\#f)(\widetilde{b}\#\widetilde{f})=\sum b_{<0>}\widetilde{b}\#((b_{<1>}%
\widetilde{f})\star f)\text{ for all }b,\widetilde{b}\in B,\text{ }f,%
\widetilde{f}\in U  \label{l-smash}
\end{equation}
(and unity $1_{B}\#\varepsilon _{H}$).
\end{punto}

\begin{definition}
Let $H$ be an $R$-bialgebra, $U\subseteq H^{\ast }$ a left $H$-module
subalgebra under the left regular $H$-action, $V\subseteq H^{\ast }$ an $R$%
-submodule and consider the $R$-linear maps 
\begin{equation}
\begin{tabular}{llllllll}
$\overline{\lambda }$ & $:$ & $H\#^{op}U$ & $\rightarrow $ & $\mathrm{End}%
_{R}(H),$ & $\sum h_{j}\otimes g_{j}$ & $\mapsto $ & $[\widetilde{k}\mapsto
\sum (g_{j}\rightharpoonup \widetilde{k})h_{j}],$ \\ 
$\overline{\rho }$ & $:$ & $V$ & $\rightarrow $ & $\mathrm{End}_{R}(H),$ & $g
$ & $\mapsto $ & $[\widetilde{k}\mapsto \widetilde{k}\leftharpoonup g].$%
\end{tabular}
\label{bar-lam-ro}
\end{equation}
We say $(V,U)$ satisfies the \emph{RL-condition with respect to }$H,$ if $%
\overline{\rho }(V)\subseteq \overline{\lambda }(H\#^{op}U),$ i.e. if 
\begin{equation}
\text{for every }g\in V,\text{ }\exists \text{ }\{(h_{j},g_{j})\}\subset
H\times U\text{ s.t. }\widetilde{k}\leftharpoonup g=\sum
(g_{j}\rightharpoonup \widetilde{k})h_{j}\text{ for every }\widetilde{k}\in
H.  \label{RL-op}
\end{equation}
\end{definition}

\begin{lemma}
\label{density-op}Let $H$ be an $R$-bialgebra, $U\subseteq H^{\ast }$ a left 
$H$-module subalgebra and consider $H$ as a right $H$-comodule algebra
through $\Delta _{H}.$ Let $\#^{op}(H,H)$ and $H\#^{op}U$ be the $R$%
-algebras defined in \emph{\ref{smash-op}} and consider the canonical $R$%
-algebra morphism $\overline{\beta }:H\#^{op}U\rightarrow \#^{op}(H,H).$

\begin{enumerate}
\item  If $_{R}H$ is finitely generated projective, then $H\#^{op}H^{\ast }%
\overset{\overline{\beta }}{\simeq }\#^{op}(H,H)$ as $R$-algebras.

\item  If $H$ is a Hopf $R$-algebra, then $\#^{op}(H,H)\simeq \mathrm{End}%
_{R}(H)^{op}$ as $R$-algebras.

\item  Let $H$ be a finitely generated projective Hopf $R$-algebra. Then $%
\overline{\lambda }:H\#^{op}H^{\ast }\rightarrow \mathrm{End}_{R}(H)^{op},$
defined in \emph{(\ref{bar-lam-ro})}, is an $R$-algebra isomorphism. In
particular $H^{\ast }$ satisfies the RL-condition \emph{(\ref{RL-op})} with
respect to $H.$

\item  If $_{R}H$ is locally projective and $U\subseteq H^{\ast }$ is \emph{%
dense}, then\emph{\ }$\overline{\beta }(H\#^{op}U)\subseteq \#^{op}(H,H)$ is
a dense $R$-subalgebra. If moreover $H$ is a projective Hopf $R$-algebra,
then $H\#^{op}U\overset{\overline{\lambda }}{\hookrightarrow }\mathrm{End}%
_{R}(H)^{op}$ is a dense $R$-subalgebra.
\end{enumerate}
\end{lemma}

\begin{Beweis}
\begin{enumerate}
\item  Since $_{R}H$ is finitely generated projective, $\overline{\beta }$
is bijective.

\item  Let $H$ be a Hopf $R$-algebra and consider the $R$-linear maps 
\begin{equation*}
\begin{tabular}{llllllllll}
$\overline{\phi }_{1}$ & $:$ & $\#^{op}(H,H)$ & $\rightarrow $ & $\mathrm{End%
}_{R}(H)^{op},$ & $f$ & $\mapsto $ & $[h$ & $\mapsto $ & $\sum
h_{1}f(h_{2})],$ \\ 
$\overline{\phi }_{2}$ & $:$ & $\mathrm{End}_{R}(H)^{op},$ & $\rightarrow $
& $\#^{op}(H,H),$ & $g$ & $\mapsto $ & $[k$ & $\mapsto $ & $\sum
S(k_{1})g(k_{2})].$%
\end{tabular}
\end{equation*}

For arbitrary $f,g\in \#^{op}(H,H)$ and $h\in H$ we have 
\begin{equation*}
\begin{tabular}{lllll}
$\overline{\phi }_{1}(f\widetilde{\star }g)(h)$ & $=$ & $\sum h_{1}(f%
\widetilde{\star }g)(h_{2})$ & $=$ & $\sum
h_{1}f(h_{3})_{1}g(h_{2}f(h_{3})_{2})$ \\ 
& $=$ & $\sum h_{11}f(h_{2})_{1}g(h_{12}f(h_{2})_{2})$ & $=$ & $\overline{%
\phi }_{1}(g)(\sum h_{1}f(h_{2}))$ \\ 
& $=$ & $(\overline{\phi }_{1}(g)\circ \overline{\phi }_{1}(f))(h),$ &  & 
\end{tabular}
\end{equation*}
i.e. $\overline{\phi }_{1}$ is an $R$-algebra morphism. For all $R$-linear
maps $f,g:H\rightarrow H$ and $h\in H$ we have 
\begin{equation*}
\begin{tabular}{lllll}
$(\overline{\phi }_{1}\circ \overline{\phi }_{2})(g)(h)$ & $=$ & $\sum h_{1}%
\overline{\phi }_{2}(g)(h_{2})$ & $=$ & $\sum h_{1}S(h_{2})g(h_{3})$ \\ 
& $=$ & $\sum \varepsilon (h_{1})g(h_{2})$ & $=$ & $g(h),$ \\ 
$(\overline{\phi }_{2}\circ \overline{\phi }_{1})(f)(h)$ & $=$ & $\sum
S(h_{1})\overline{\phi }_{1}(f)(h_{2})$ & $=$ & $\sum S(h_{1})h_{2}f(h_{3})$
\\ 
& $=$ & $\sum \varepsilon (h_{1})f(h_{2})$ & $=$ & $f(h).$%
\end{tabular}
\end{equation*}
\newline
Hence $\overline{\phi }_{1}$ is an $R$-algebra isomorphism with inverse $%
\overline{\phi }_{2}.$

\item  Let $H$ be a finitely generated projective Hopf $R$-algebra. By (1)
and (2) $H\#^{op}H^{\ast }\overset{\overline{\beta }}{\simeq }\#^{op}(H,H)%
\overset{\overline{\phi }_{1}}{\simeq }\mathrm{End}_{R}(H)^{op}$ as $R$%
-algebras. Hence $\overline{\lambda }=\overline{\phi }_{1}\circ \overline{%
\beta }:H\#^{op}H^{\ast }\rightarrow \mathrm{End}_{R}(H)^{op}$ is an $R$%
-algebra isomorphism. In particular $\overline{\rho }(H^{\ast })\subseteq 
\mathrm{End}_{R}(H)^{op}=\overline{\lambda }(H\#^{op}H^{\ast }),$ i.e. $%
H^{\ast }$ satisfies the RL-condition (\ref{RL-op}) with respect to $H.$

\item  By \cite[Theorem 3.18 (2)]{Abu2002} $\overline{\beta }%
(H\#^{op}U)\subseteq \#^{op}(H,H)$ is a dense $R$-subalgebra. If $H$ is a
Hopf $R$-algebra, then $\#(H,H)\overset{\overline{\phi }_{1}}{\simeq }%
\mathrm{End}_{R}(H)^{op}$ as $R$-algebras by (2) and we are done (notice
that $\overline{\beta }$ is an embedding, if $_{R}H$ is projective).$%
\blacksquare $
\end{enumerate}
\end{Beweis}

\begin{lemma}
\label{alp=bet-op}Let $H$ be a Hopf $R$-algebra, $A$ an $R$-algebra, $%
U\subseteq H^{\ast }$ a left $H$-submodule and consider the $R$-paring $%
\overline{P}:=(H,U).$ Then the canonical $R$-linear map $\overline{\alpha }%
:=\alpha _{A\otimes _{R}H}^{\overline{P}}:(A\otimes _{R}H)\otimes
_{R}U\rightarrow \mathrm{Hom}_{R}(H,A\otimes _{R}H)$ is injective if and
only if the following map is injective 
\begin{equation}
\overline{\chi }:A\otimes _{R}(H\otimes _{R}U)\rightarrow \mathrm{End}%
_{A-}(A\otimes _{R}H),\text{ }a\otimes (h\otimes f)\mapsto \lbrack (%
\widetilde{a}\otimes k)\mapsto \widetilde{a}a\otimes (f\rightharpoonup k)h].
\label{bet-bar}
\end{equation}
\end{lemma}

\begin{Beweis}
First we show that the $R$-linear map 
\begin{equation*}
\overline{\epsilon }:\mathrm{Hom}{\normalsize _{R}}(H,A\otimes
_{R}H)\rightarrow \mathrm{End}_{A-}(A\otimes _{R}H),\text{ }g\mapsto \lbrack 
\widetilde{a}\otimes k\mapsto (\widetilde{a}\otimes k_{1})g(k_{2})]
\end{equation*}
is bijective with inverse 
\begin{equation*}
\overline{\epsilon }^{-1}:\mathrm{End}_{A-}(A\otimes _{R}H)\rightarrow 
\mathrm{Hom}{\normalsize _{R}}(H,A\otimes _{R}H),\text{ }f\mapsto \lbrack
k\mapsto (1_{A}\otimes S(k_{1}))f(1_{A}\otimes k_{2})].
\end{equation*}
In fact we have for all $f\in \mathrm{End}_{A-}(A\otimes _{R}H),$ $k\in H,$ $%
\widetilde{a}\in A:$%
\begin{equation*}
\begin{tabular}{lll}
$\overline{\epsilon }(\overline{\epsilon }^{-1}(f))(\widetilde{a}\otimes k)$
& $=$ & $\sum (\widetilde{a}\otimes k_{1})\overline{\epsilon }^{-1}(f)(k_{2})
$ \\ 
& $=$ & $\sum (\widetilde{a}\otimes k_{1})(1_{A}\otimes
S(k_{2}))f(1_{A}\otimes k_{3})$ \\ 
& $=$ & $\sum (\widetilde{a}\otimes k_{1}S(k_{2}))f(1_{A}\otimes k_{3})$ \\ 
& $=$ & $\sum (\widetilde{a}\otimes \varepsilon
_{H}(k_{1})1_{H})f(1_{A}\otimes k_{2})$ \\ 
& $=$ & $\sum (\widetilde{a}\otimes 1_{H})f(1_{A}\otimes k)$ \\ 
& $=$ & $f(\widetilde{a}\otimes k)$%
\end{tabular}
\end{equation*}
and for all $g\in \mathrm{Hom}_{R}(H,A\otimes _{R}H),$ $k\in H:$%
\begin{equation*}
\begin{tabular}{lll}
$\overline{\epsilon }^{-1}(\overline{\epsilon }(g))(k)$ & $=$ & $\sum
(1_{A}\otimes S(k_{1}))\overline{\epsilon }(g)(1_{A}\otimes k_{2})$ \\ 
& $=$ & $\sum (1_{A}\otimes S(k_{1}))(1\otimes k_{2})g(k_{3})$ \\ 
& $=$ & $\sum (1_{A}\otimes S(k_{1})k_{2})g(k_{3})$ \\ 
& $=$ & $\sum (1_{A}\otimes \varepsilon _{H}(k_{1})1_{H})g(k_{2})$ \\ 
& $=$ & $(1_{A}\otimes 1_{H})g(k)$ \\ 
& $=$ & $g(k).$%
\end{tabular}
\end{equation*}

Moreover we have for all $a\in A,$ $h\in H,$ $f\in U$ and $k\in H:$%
\begin{equation*}
\begin{tabular}{lll}
$(\overline{\epsilon }\circ \overline{\alpha })(a\otimes (h\otimes f))(%
\widetilde{a}\otimes k)$ & $=$ & $\sum (\widetilde{a}\otimes k_{1})\alpha
_{A\otimes _{R}H}^{\overline{P}}(a\otimes (h\otimes f))(k_{2}))$ \\ 
& $=$ & $\sum (\widetilde{a}\otimes k_{1})(a\otimes h)f(k_{2})$ \\ 
& $=$ & $\sum \widetilde{a}a\otimes k_{1}f(k_{2})h$ \\ 
& $=$ & $\widetilde{a}a\otimes (f\rightharpoonup k)h$ \\ 
& $=$ & $\overline{\chi }(a\otimes (h\otimes f))(\widetilde{a}\otimes k),$%
\end{tabular}
\end{equation*}
i.e. $\overline{\chi }=\overline{\epsilon }\circ \overline{\alpha }.$
Consequently $\overline{\chi }$ is injective iff $\overline{\alpha }$ is so.$%
\blacksquare $
\end{Beweis}

\begin{punto}
Let $H$ be a Hopf $R$-algebra, $A\#_{\sigma }H$ a right $H$-crossed product
with invertible cocycle and consider the $R$-linear maps $\overline{\varphi }%
,\overline{\psi }:H\otimes _{R}A\rightarrow \mathrm{Hom}_{R}(H,A)$ defined
as 
\begin{equation*}
\begin{tabular}{lll}
$\overline{\varphi }(h\otimes a)(\widetilde{h})$ & $=$ & $\sum [\widetilde{h}%
_{1}a]\sigma (\widetilde{h}_{2}\otimes h),$ \\ 
$\overline{\psi }(h\otimes a)(\widetilde{h})$ & $=$ & $\sum \sigma ^{-1}(S(%
\widetilde{h}_{3})\otimes \widetilde{h}_{4})[S(\widetilde{h}_{2})a]\sigma (S(%
\widetilde{h}_{1})\otimes \widetilde{h}_{5}h).$%
\end{tabular}
\end{equation*}
Let $U\subseteq H^{\ast }$ be a \emph{left }$H$-module subalgebra, $%
V\subseteq H^{\ast }$ an $R$-submodule and consider the $R$-linear map $%
J:A\otimes _{R}V\rightarrow \mathrm{Hom}_{R}(H,A).$ We say $(V,U)$ is \emph{%
compatible}, if the following conditions are satisfied:

\begin{enumerate}
\item  $\overline{\varphi }(H\otimes _{R}A),$ $\overline{\psi }(H\otimes
_{R}A)\subseteq J(A\otimes _{R}V);$

\item  $(V,U)$ satisfies the \emph{RL-condition} (\ref{RL-op}) with respect
to $H.$
\end{enumerate}
\end{punto}

Analog to Proposition \ref{Chen} and in the light of Lemma \ref{alp=bet-op}
and the modified RL-condition (\ref{RL-op}) we restate \cite[Theorem 8,
Corollary 9]{Che93} for the case of an arbitrary commutative ground ring:

\begin{proposition}
\label{Chen-op}Let $H$ be a Hopf $R$-algebra, $A\#_{\sigma }H$ be a right $H$%
-crossed product with invertible cocycle, $U\subseteq H^{\ast }$ a left $H$%
-module subalgebra and consider the $R$-pairing $\overline{P}:=(H,U).$
Assume there exists an $R$\emph{-submodule} $V\subseteq H^{\ast },$ such
that $(V,U)$ is compatible. If the canonical $R$-linear map $\overline{%
\alpha }:=\alpha _{A\otimes _{R}H}^{\overline{P}}:(A\otimes _{R}H)\otimes
_{R}U\rightarrow \mathrm{Hom}_{R}(H,A\otimes _{R}H)$ is injective, then
there exists an $R$-algebra isomorphism 
\begin{equation*}
(A\#_{\sigma }H)\#^{op}U\simeq A\otimes _{R}(H\#^{op}U).
\end{equation*}
\end{proposition}

\begin{Beweis}
By \cite[Lemma 7]{Che93} we have a commutative diagram of $R$-algebra
morphisms 
\begin{equation}
\xymatrix{ & (A \#_{\sigma} H) \# ^{op} U \ar[dl]_{\overline{\alpha }}
\ar[dr]^{\overline {\gamma} } & \\ \# ^{op} (H,A \#_{\sigma} H)
\ar[rr]^{\overline{\pi }} & & {\rm End}_{A-} (A \otimes_{R} H)^{op} \\ & A
\otimes_{R} (H \# ^{op} U) \ar[ur]_{\overline {\chi}} \ar[ul]^{\overline
{\delta } } & }  \label{alg-d-op}
\end{equation}
where 
\begin{equation*}
\begin{tabular}{lll}
$\overline{\alpha }((a\#h)\#^{op}f)(k)$ & $=$ & $(a\#h)f(k),$ \\ 
$\overline{\chi }(a\otimes (h\#^{op}f))(\widetilde{a}\otimes k)$ & $=$ & $%
\widetilde{a}a\otimes (f\rightharpoonup k)h,$ \\ 
$\overline{\gamma }((a\#h)\#^{op}f)(\widetilde{a}\otimes k)$ & $=$ & $\sum 
\widetilde{a}[k_{1}a]\sigma (k_{2}\otimes h_{1})\otimes (f\rightharpoonup
k_{3})h_{2},$ \\ 
$\overline{\delta }(a\otimes (h\#^{op}f))(k)$ & $=$ & $\sum \sigma
^{-1}(S(k_{4})\otimes k_{5})[S(k_{3})a]\sigma (S(k_{2})\otimes
k_{6}h_{1})\#S(k_{1})(f\rightharpoonup k_{7})h_{2},$ \\ 
$\overline{\pi }(g)(\widetilde{a}\otimes k)$ & $=$ & $\sum (\widetilde{a}%
\#_{\sigma }k_{1})g(k_{2}).$%
\end{tabular}
\end{equation*}
By assumption $\overline{\alpha }:=\alpha _{A\otimes _{R}H}^{\overline{P}}$
is injective, hence $\overline{\chi }$ is by Lemma \ref{alp=bet-op}
injective and consequently $\overline{\delta }$ is injective. Moreover $%
\overline{\pi }$ is an $R$-algebra isomorphism by \cite[Propsoition 4.1]
{Kop92}, hence $\overline{\gamma }$ is injective. It remains then to show
that $\mathrm{\func{Im}}(\overline{\gamma })\subseteq \mathrm{\func{Im}}(%
\overline{\chi })$ and $\mathrm{\func{Im}}(\overline{\delta })\subseteq 
\mathrm{\func{Im}}(\overline{\alpha }).$ For arbitrary $a\otimes h\in
A\otimes _{R}H,$ there exists $\sum a_{u}\otimes g_{u}\in A\otimes _{R}V$
such that $\overline{\varphi }(h_{1}\otimes a)=J(\sum a_{u}\otimes g_{u})$
and moreover there exists $\sum h_{uj}\#^{op}g_{uj}\in H\otimes _{R}U$ with $%
\overline{\rho }(g_{u})=\overline{\lambda }(\sum h_{uj}\#^{op}g_{uj}).$ So
for all $a,\widetilde{a}\in A,$ $h,k\in H$ and $f\in U:$%
\begin{equation*}
\begin{tabular}{lllll}
$\overline{\gamma }((a\#h)\#^{op}f)(\widetilde{a}\otimes k)$ & $=$ & $\sum 
\widetilde{a}[k_{1}a]\sigma (k_{2}\otimes h_{1})\otimes (f\rightharpoonup
k_{3})h_{2}$ &  &  \\ 
& $=$ & $\sum \widetilde{a}[k_{11}a]\sigma (k_{12}\otimes h_{1})\otimes
(f\rightharpoonup k_{2})h_{2}$ &  &  \\ 
& $=$ & $\sum \widetilde{a}\overline{\varphi }(h_{1}\otimes a)(k_{1})\otimes
(f\rightharpoonup k_{2})h_{2}$ &  &  \\ 
& $=$ & $\sum \widetilde{a}J(\sum a_{u}\otimes g_{u})(k_{1})\otimes
(f\rightharpoonup k_{2})h_{2}$ &  &  \\ 
& $=$ & $\sum \widetilde{a}a_{u}g_{u}(k_{1})\otimes (f\rightharpoonup
k_{2})h_{2}$ &  &  \\ 
& $=$ & $\sum \widetilde{a}a_{u}\otimes g_{u}(k_{1})k_{2}f(k_{3})h_{2}$ &  & 
\\ 
& $=$ & $\sum \widetilde{a}a_{u}\otimes (k_{1}\leftharpoonup
g_{u})f(k_{2})h_{2}$ &  &  \\ 
& $=$ & $\sum \widetilde{a}a_{u}\otimes (g_{u,j}\rightharpoonup
k_{1})h_{u,j}f(k_{2})h_{2}$ &  &  \\ 
& $=$ & $\sum \widetilde{a}a_{u}\otimes
k_{1}g_{u,j}(k_{2})f(k_{3})h_{u,j}h_{2}$ &  &  \\ 
& $=$ & $\sum \widetilde{a}a_{u}\otimes ((g_{u,j}\star f)\rightharpoonup
k)h_{u,j}h_{2}$ &  &  \\ 
& $=$ & $\overline{\chi }(a_{u}\otimes (h_{u,j}h_{2}\#^{op}(g_{u,j}\star
f)))(\widetilde{a}\otimes k),$ &  & 
\end{tabular}
\end{equation*}
i.e. $\mathrm{\func{Im}}(\overline{\gamma })\subseteq \mathrm{\func{Im}}(%
\overline{\chi }).$ For arbitrary $a\otimes h\in A\otimes _{R}H,$ there
exists $\sum a_{w}\otimes g_{w}\in A\otimes _{R}V$ such that $\overline{%
\varphi }(h_{1}\otimes a)=J(\sum a_{w}\otimes g_{w})$ and moreover there
exists $\sum h_{wj}\#^{op}g_{wj}\in H\otimes _{R}U$ with $\overline{\rho }%
(g_{u})=\overline{\lambda }(\sum h_{wj}\#^{op}g_{wj}).$ So we have for all $%
a\in A,$ $h,k\in H$ and $f\in U:$ 
\begin{equation*}
\begin{tabular}{lll}
$\overline{\delta }(a\otimes (h\#^{op}f))(k)$ & $=$ & $\sum \sigma
^{-1}(S(k_{4})\otimes k_{5})[S(k_{3})a]\sigma (S(k_{2})\otimes
k_{6}h_{1})\#S(k_{1})(f\rightharpoonup k_{7})h_{2}$ \\ 
& $=$ & $\sum \sigma ^{-1}(S(k_{23})\otimes k_{24})[S(k_{22})a]\sigma
(S(k_{21})\otimes k_{25}h_{1})\#S(k_{1})(f\rightharpoonup k_{3})h_{2}$ \\ 
& $=$ & $\sum \overline{\psi }(h_{1}\otimes
a)(k_{2})\#S(k_{1})(f\rightharpoonup k_{3})h_{2}$ \\ 
& $=$ & $\sum J(\sum a_{w}\otimes g_{w})(k_{2})\#S(k_{1})(f\rightharpoonup
k_{3})h_{2}$ \\ 
& $=$ & $\sum a_{w}g_{w}(k_{2})\#S(k_{1})(f\rightharpoonup k_{3})h_{2}$ \\ 
& $=$ & $\sum a_{w}\#S(k_{1})g_{w}(k_{2})k_{3}f(k_{4})h_{2}$ \\ 
& $=$ & $\sum a_{w}\#S(k_{1})(k_{2}\leftharpoonup g_{w})f(k_{3})h_{2}$ \\ 
& $=$ & $\sum a_{w}\#S(k_{1})(g_{w,j}\rightharpoonup
k_{2})h_{w,j}f(k_{3})h_{2}$ \\ 
& $=$ & $\sum a_{w}\#S(k_{1})k_{2}g_{w,j}(k_{3})f(k_{4})h_{w,j}h_{2}$ \\ 
& $=$ & $\sum a_{w}\#g_{w,j}(k_{1})f(k_{2})h_{w,j}h_{2}$ \\ 
& $=$ & $\sum a_{w}\#(g_{w,j}\star f)(k)h_{w,j}h_{2}$ \\ 
& $=$ & $\overline{\alpha }((a_{w}\#h_{w,j}h_{2})\#^{op}g_{w,j}\star f)(k),$%
\end{tabular}
\end{equation*}
i.e. $\mathrm{\func{Im}}(\overline{\delta })\subseteq \mathrm{\func{Im}}(%
\overline{\alpha })$ and we are done.$\blacksquare $
\end{Beweis}

As a consequence of Lemma \ref{M-pure} and Proposition \ref{Chen-op} we get
an analog to Theorem \ref{haupt}, which generalizes \cite[Corollary 9]{Che93}
(resp. \cite[Corollary 10]{Che93}) from the case of a base field (resp. a
Dedekind domain) to the case of an arbitrary Noetherian ground ring:

\begin{theorem}
\label{haupt-op}Let $R$ be Noetherian, $H$ a Hopf $R$-algebra, $A\#_{\sigma
}H$ a right $H$-crossed product with invertible cocycle, $U\subseteq H^{\ast
}$ a left $H$-module subalgebra and consider the $R$-pairing $\overline{P}%
:=(H,U).$ Assume there exists an $R$-submodule $V\subseteq H^{\ast },$ such
that $(V,U)$ is compatible. If $U\subset R^{H}$ is $A\otimes _{R}H$-pure 
\emph{(}e.g. $H$ is a Hopf $\alpha $-algebra and $U\subseteq H^{\circ }$ is
an $R$-subbialgebra\emph{)}, then we have an $R$-algebra isomorphism 
\begin{equation*}
(A\#_{\sigma }H)^{op}\#U\simeq A\otimes _{R}(H\#^{op}U).
\end{equation*}
\end{theorem}

\begin{corollary}
\label{dense-op}Let $H$ be a projective Hopf $R$-algebra, $A\#_{\sigma }H$ a
right $H$-crossed product with invertible cocycle, $U\subseteq H^{\ast }$ a
left $H$-module subalgebra and consider the $R$-paring $P:=(H,U).$ Assume
there exists an $R$-submodule $V\subseteq H^{\ast },$ such that $(V,U)$ is
compatible. If $U\subseteq H^{\ast }$ is dense and the canonical $R$-linear
map $\alpha _{A\otimes _{R}H}^{P}$ is injective \emph{(}e.g. $R$ is
Noetherian and $U\subseteq R^{H}$ is $A$-pure\emph{)}, then there is a dense%
\emph{\ }$R$-subalgebra $\mathcal{L}\subseteq \mathrm{End}_{R}(H)^{op}$ and
an $R$-algebra isomorphism 
\begin{equation*}
(A\#_{\sigma }H)\#^{op}U\simeq A\otimes _{R}\mathcal{L}.
\end{equation*}
This is the case in particular, if $R$ is a QF ring, $H$ is a residually
finite Hopf $\alpha $-algebra and $U\subseteq H^{\circ }$ is a dense $R$%
-subbialgebra.
\end{corollary}

\begin{Beweis}
If $U\subseteq H^{\ast }$ is dense, then $\mathcal{L}:=H\#^{op}U\overset{%
\overline{\lambda }}{\hookrightarrow }\mathrm{End}_{R}(H)^{op}$ is by Lemma 
\ref{density-op} a dense $R$-subalgebra. If $\alpha _{A\otimes _{R}H}^{P}$
is injective, then the isomorphism follows by Theorem \ref{haupt-op}. If $R$
is a QF ring and $H$ is a residually finite Hopf $\alpha $-algebra, then $%
H^{\circ }\subset H^{\ast }$ is dense by \cite[Proposition 2.4.19]{Abu2001}.
If moreover $U\subseteq H^{\circ }$ is a dense $R$-subbialgebra, then $%
U\subseteq H^{\ast }$ is dense, $\alpha _{A\otimes _{R}H}^{P}$ is injective
and we are done.$\blacksquare $
\end{Beweis}

Similar argument to those in the proof of Corollary \ref{f.d.} can be used
to prove

\begin{corollary}
\label{f.d.-op}Let $H$ be a Hopf $R$-algebra and $A\#_{\sigma }H$ a right $H$%
-crossed product with invertible cocycle. Then we have an isomorphism of $R$%
-algebras 
\begin{equation*}
(A\#_{\sigma }H)\#^{op}H^{\ast }\simeq A\otimes _{R}(H\#^{op}H^{\ast })
\end{equation*}
at least when

\begin{enumerate}
\item  $_{R}H$ is finitely generated projective, \emph{or}

\item  $_{R}A$ is finitely generated, $H$ is cocommutative and $\alpha
_{A\otimes _{R}H}^{P}$ is injective \emph{(}e.g. $R$ is Noetherian and $%
H^{\ast }\hookrightarrow R^{H}$ is $A\otimes _{R}H$-pure\emph{)}.
\end{enumerate}
\end{corollary}

\section*{The subalgebra $H^{\protect\omega }\subseteq H^{\ast }$}

In what follows let $H$ be a \emph{locally projective} Hopf $R$-algebra and
consider the measuring $\alpha $-pairing $P:=(H^{\ast },H)$ (notice that the
canonical $R$-linear map $\alpha _{R}^{P}:H\rightarrow H^{\ast \ast }$ is
injective).

\begin{punto}
Consider $H^{\ast }$ with the right $H^{\ast }$-action 
\begin{equation*}
(f\leftharpoonup g)(h):=\sum f(h_{2})g(S(h_{1})h_{3})\text{ for all }f,g\in
H^{\ast }\text{ and }h\in H.
\end{equation*}
Then $H^{\ast }$ is a right $H^{\ast }$-module and $H^{\omega }:=$ $^{H}%
\mathrm{Rat}(H_{H^{\ast }}^{\ast })$ is analog to Theorem \ref{equal} a left 
$H$-comodule with structure map $\omega :H^{\omega }\rightarrow H\otimes
_{R}H^{\omega }.$
\end{punto}

Analog to \cite[Propositions 3.2, 3.3]{Kop92} we have

\begin{proposition}
\label{wH-prop}Consider the left $H$-comodule $(H^{\omega },\omega ).$

\begin{enumerate}
\item  If $f\in H^{\omega },$ then $\omega (f)=\sum f_{<-1>}\otimes f_{<0>}$
satisfies the following conditions:

\begin{enumerate}
\item  $f\star g=\sum f_{<-1>}g\star f_{<0>}$ for all $g\in H^{\ast }.$

\item  $h\leftharpoonup f=\sum (f_{<0>}\rightharpoonup h)f_{<-1>}$ for all $%
h\in H.$

\item  $\sum f(h_{2})S(h_{1})h_{3}=\sum f_{<0>}(h)f_{<-1>}$ for all $h\in H.$
\end{enumerate}

\item  Let $f\in H^{\ast }.$ If there exists $\zeta =\sum f_{<-1>}\otimes
f_{<0>}\in H\otimes _{R}H^{\ast }$ that satisfies any of the conditions in 
\emph{(1)}, then $f\in H^{\omega }$ and $\omega (f)=\zeta .$

\item  $H^{\omega }\subseteq H^{\ast }$ is an $R$-subalgebra and moreover a
left $H$-comodule algebra.

\item  $H^{\omega }\subseteq H^{\ast }$ is a left $H$-module subalgebra with 
\begin{equation*}
\omega (hf)=\sum h_{1}f_{<-1>}S(h_{3})\otimes h_{2}f_{<0>}\text{ for all }%
h\in H\text{ and }f\in H^{\omega }.
\end{equation*}
\end{enumerate}
\end{proposition}

\qquad As a consequence of Proposition \ref{Chen-op} and Theorem \ref
{haupt-op} we get the following generalization of \cite[Theorem 4.2]{Kop92}:

\begin{theorem}
\label{ow-1}Let $H$ be a locally projective Hopf $R$-algebra, $U\subseteq
H^{\ast }$ a left $H$-module subalgebra, $P:=(H,U)$ the induced $R$-pairing
and assume that $\overline{\varphi }(H\otimes _{R}A),$ $\overline{\psi }%
(H\otimes _{R}A)\subseteq J(A\otimes _{R}\omega ^{-1}(H\otimes _{R}U)).$ If $%
\alpha _{A\otimes _{R}H}^{P}$ is injective \emph{(}e.g. $R$ is Noetherian
and $U\subset R^{H}$ is $A$-pure\emph{)}, then there is an $R$-algebra
isomorphism 
\begin{equation*}
(A\#_{\sigma }H)^{op}\#U\simeq A\otimes _{R}(H\#^{op}U).
\end{equation*}
\end{theorem}

\begin{Beweis}
Consider the $R$-submodule $V:=\omega ^{-1}(H\otimes _{R}U).$ Since $%
V\subseteq H^{\omega },$ it's clear by Proposition \ref{wH-prop} (1-b) that $%
(V,U)$ satisfies the RL-condition (\ref{RL-op}) with respect to $H.$
Consequently $(V,U)$ is compatible. If $\alpha _{A\otimes _{R}H}^{P}$ is
injective, then we are done by Proposition \ref{Chen-op}.$\blacksquare $
\end{Beweis}

\begin{corollary}
\label{normal}Let $H$ be a locally projective Hopf $R$-algebra, $U\subseteq
H^{\omega }$ a left $H$-module subalgebra, $P:=(H,U)$ the induced $R$%
-pairing and assume that $\omega (U)\subseteq H\otimes _{R}U$ and $\overline{%
\varphi }(H\otimes _{R}A),$ $\overline{\psi }(H\otimes _{R}A)\subseteq
J(A\otimes _{R}U).$ If $\alpha _{A\otimes _{R}H}^{P}$ is injective \emph{(}%
e.g. $R$ is Noetherian and $U\subset R^{H}$ is $A$-pure\emph{)}, then there
is an $R$-algebra isomorphism 
\begin{equation*}
(A\#_{\sigma }H)^{op}\#U\simeq A\otimes _{R}(H\#^{op}U).
\end{equation*}
\end{corollary}

\begin{remark}
If the Hopf algebra $H$ has a bijective antipode then it has a twisted
antipode, namely $\overline{S}:=S^{-1}.$ In the proofs (by different
authors) of several duality theorems for smash products assuming the
bijectivity of the antipode, no use was made of $S\circ
S^{-1}=id=S^{-1}\circ S;$ instead there was a heavy use of the main
properties of $S^{-1},$ namely that it is an algebra and coalgebra
anti-morphism, and mainly that
\begin{equation*}
\sum S^{-1}(h_{2})h_{1}=\varepsilon (h)1_{H}=\sum h_{2}S^{-1}(h_{1})\text{
for every }h\in H.
\end{equation*}
A twisted antipode has also these main properties and this is why the
original versions (in \cite{Abu2001}) of the results in section two remain
true after replacing the bijectivity of the antipode by the weaker condition
of the existence of a twisted antipode!!
\end{remark}

\begin{punto}
\label{op-rel}(Compare \cite[Lemma 5.3]{Kop95})\emph{\ }Let $H$ be a Hopf $R$%
-algebra with a twisted antipode $\overline{S}$ and $A\#_{\sigma }H$ a right 
$H$-crossed product with invertible cocycle $\sigma .$ Then $h^{op}a^{op}:=%
\overline{S}(h)a$ induces on $A^{op}$ a weak left $H^{op}$-action and $%
A^{op}\#_{\tau }H^{op}$ is a right $H^{op}$-crossed product with invertible
cocycle 
\begin{equation*}
\tau :H\otimes _{R}H\rightarrow A,\text{ }(h,k)\mapsto \sigma ^{-1}(%
\overline{S}(h),\overline{S}(k)).
\end{equation*}
Moreover $A\#_{\sigma }H\simeq (A^{op}\#_{\tau }H^{op})^{op}$ as right $H$%
-comodule algebras.
\end{punto}

\begin{remark}
As indicated earlier, the original versions (\cite{Abu2001}) of the main
duality theorems for smash products were proved under the assumption of the
bijectivity of the antipode of $H$ and it was not clear why such an
assumption is not needed in the corresponding results for opposite smash
products. Upon suggestion of the referee this condition is replaced in this
paper with the weaker condition that $H$ has a twisted antipode which
clarifies, to some extent, this issue (notice that the rule of $H$ is played
in the third section by $H^{op}$ which has a twisted antipode!!). However,
it should be noted that the results in the third section cannot be deduced
directly from the corresponding results in the second section, since (in
light of \ref{op-rel}) we have to assume that $H$ has a twisted antipode!!
\end{remark}

\qquad However, some of duality theorems for smash products can be deduced
from the corresponding ones for opposite smash products under the assumption
that $H$ has a twisted antipode, for example we have

\begin{punto}
Let $R$ be Noetherian, $H$ a Hopf $\alpha $-algebra with twisted antipode, $%
U\subseteq H^{\circ }$ an $R$-subbialgebra and consider the $R$-subbialgebra 
$U^{cop}\subseteq \left( H^{op}\right) ^{\circ }.$ Assume there exists an $R$%
-submodule $V\subseteq (H^{op})^{\ast },$ such that 
\begin{equation}
\text{for every }g\in V,\text{ there exist }\{(h_{j},g_{j})\}\subset H\times
U,\text{ s.t. }\widetilde{h}\leftharpoonup g=\sum h_{j}(g_{j}\rightharpoonup 
\widetilde{h})\text{ for all }h\in H\text{ }
\end{equation}
and that for every $(h,a)\in H\times A$ there exist subclasses $%
\{a_{u},g_{u}\},\{b_{w},g_{w}\}\subset A\times V$ with 
\begin{equation*}
\begin{tabular}{lll}
$\sum \sigma ^{-1}(\overline{S}(\widetilde{h}_{2})\otimes \overline{S}(h))[%
\overline{S}(\widetilde{h}_{1})a]$ & $=$ & $\sum a_{u}g_{u}(\widetilde{h}),$
\\ 
$\sum \sigma ^{-1}(\widetilde{h}_{1},\overline{S}(\widetilde{h}_{5}h))[%
\widetilde{h}_{2}a]\sigma (\widetilde{h}_{3}\otimes \overline{S}(\widetilde{h%
}_{4}))$ & $=$ & $\sum b_{w}g_{w}(\widetilde{h}).$%
\end{tabular}
\end{equation*}
Combining \cite[Corollary 2.4]{Kop92} and Theorem \ref{haupt-op} we get the $%
R$-algebra isomorphisms 
\begin{equation*}
\begin{tabular}{lllll}
$(A\#_{\sigma }H)\#U$ & $\simeq $ & $((A^{op}\#_{\tau
}H^{op})\#^{op}U^{cop})^{op}$ & $\simeq $ & $A\otimes
_{R}(H^{op}\#^{op}U^{cop})^{op}$ \\ 
& $\simeq $ & $(A^{op}\otimes _{R}(H^{op}\#^{op}U^{cop}))^{op}$ & $\simeq $
& $A\otimes _{R}(H\#U).$%
\end{tabular}
\end{equation*}
\end{punto}

\bigskip \qquad 

\textbf{Acknowledgments.} Part of this note \emph{extends} and \emph{improves%
} some results of my doctoral thesis at the Heinrich-Heine Universit\"{a}t, D%
\"{u}sseldorf (Germany). I am so grateful to my advisor Prof. Robert
Wisbauer for his wonderful supervision and the continuous encouragement and
support. I am also thankful to the referee for his/her useful suggestions.
Many thanks to KFUPM for the financial support and the excellent research
facilities.

\end{document}